\pdfoutput=1
\documentclass[11pt,twoside]{article}
\usepackage[letterpaper,margin=0.9in]{geometry}
\usepackage{setspace}
\usepackage{amsmath,amssymb,amsthm,mathtools}
\usepackage{graphicx,booktabs,enumitem}
\usepackage{algorithm,algpseudocode}
\usepackage{titlesec,fancyhdr,microtype,placeins,array}
\usepackage[dvipsnames]{xcolor}
\usepackage[colorlinks=true,linkcolor=blue!50!black,citecolor=green!40!black,urlcolor=blue!50!black,
  pdftitle={Supermodularity and Robust Optimization in Elliptic PDEs},
  pdfauthor={Debangana Mukherjee and Akash Ashirbad Panda}]{hyperref}
\usepackage[capitalise]{cleveref}
\usepackage{orcidlink}

\usepackage{titlesec}

\titleformat{\section}[block]
  {\normalfont\large\bfseries\centering}
  {\thesection.}
  {0.6em}
  {}

\titlespacing*{\section}
  {0pt}
  {2.5ex plus 0.6ex minus 0.2ex}
  {1.2ex}

\titleformat{\subsection}[block]
  {\normalfont\large\bfseries}
  {\thesubsection.}
  {0.5em}
  {}

\titlespacing*{\subsection}
  {0pt}
  {2.0ex plus 0.5ex minus 0.2ex}
  {0.8ex}

\titleformat{\paragraph}[runin]
  {\normalfont\itshape}
  {}
  {0em}
  {}
  [.]

\titlespacing*{\paragraph}
  {0pt}
  {1.2ex plus 0.3ex}
  {0.5em}

\theoremstyle{plain}
\newtheorem{theorem}{Theorem}[section]
\newtheorem{lemma}[theorem]{Lemma}
\newtheorem{proposition}[theorem]{Proposition}
\newtheorem{corollary}[theorem]{Corollary}
\theoremstyle{definition}
\newtheorem{definition}[theorem]{Definition}
\newtheorem{assumption}[theorem]{Assumption}
\newtheorem{example}[theorem]{Example}
\theoremstyle{remark}
\newtheorem{remark}[theorem]{Remark}
\numberwithin{equation}{section}
\crefname{assumption}{Assumption}{Assumptions}
\crefname{figure}{Figure}{Figures}
\crefname{table}{Table}{Tables}
\crefname{algorithm}{Algorithm}{Algorithms}
\crefname{appendix}{Appendix}{Appendices}

\renewenvironment{abstract}
{%
  \begin{center}
  \begin{minipage}{0.88\textwidth}
  \small
  \noindent\textbf{Abstract.}\ \ignorespaces
}
{%
  \end{minipage}
  \end{center}
  \smallskip
}
\newenvironment{keywords}{\par\small\noindent\textbf{Key words.}\ \ignorespaces}{\par\smallskip}
\newenvironment{MSCcodes}{\par\small\noindent\textbf{MSC codes.}\ \ignorespaces}{\par\medskip}

\usepackage{etoolbox}

\AtBeginEnvironment{thebibliography}{%
  \small
  \setlength{\itemsep}{1pt}%
  \setlength{\parskip}{0pt}%
}

\makeatletter
\def\@fnsymbol#1{\ensuremath{\ifcase#1\or *\or \dagger\or \ddagger\or \mathsection\or \mathparagraph\or \|\or **\or \dagger\dagger\or \ddagger\ddagger\else\@ctrerr\fi}}
\makeatother

\newcommand{\R}{\mathbb{R}}
\newcommand{\E}{\mathcal{E}}
\newcommand{\Pcal}{\mathcal{P}}
\newcommand{\dmu}{\,\mathrm{d}\mu}
\newcommand{\dx}{\,\mathrm{d}x}
\newcommand{\Vp}{V_{+}}
\newcommand{\Vsp}{V^{*}_{+}}
\newcommand{\xih}{\widehat\xi}
\newcommand{\xib}{\bar\xi}
\newcommand{\xil}{\underline\xi}
\newcommand{\ip}[2]{\langle #1,#2\rangle}
\newcommand{\Smax}{S^{\max}}

\title{\large\bfseries SUPERMODULARITY AND ROBUST OPTIMIZATION IN ELLIPTIC PDES}
\author{\normalsize\scshape Debangana Mukherjee\orcidlink{0000-0002-1516-840X}\thanks{School of Interwoven Arts and Sciences, Krea University, Sri City,
India (\href{mailto:debangana.mukherjee@krea.edu.in}{\nolinkurl{debangana.mukherjee@krea.edu.in}}).}
\and \normalsize\scshape Akash Ashirbad Panda\orcidlink{0000-0003-3496-9111}\thanks{Corresponding author. Department of Mathematics, Indian Institute of
Technology Indore, Simrol, India (\href{mailto:akashpanda@iiti.ac.in}{\nolinkurl{akashpanda@iiti.ac.in}}).}}
\date{}

\begin{document}
\maketitle
\thispagestyle{plain}

\begin{abstract}
We study distributionally robust optimization problems constrained by semilinear elliptic PDEs. Uncertain
nonnegative loads enter the PDE through $m$ intensity parameters, and damage is measured by integrating a
nondecreasing convex function of the solution. We prove that, for coercive operators satisfying a
positivity-preserving lattice condition, a nondecreasing concave reaction $g$ makes the damage nondecreasing,
convex, and supermodular in the load intensities. Under nonnegative loading and $g(0)=0$, concavity is also
necessary for these properties to hold for all admissible operators and damage functions. For Wasserstein
uncertainty sets based on a weighted $\ell_1$ distance and $N$ observations, we show that at most $m+1$ nested
stress corners are needed per observation. These corners form a set of at most $N(m+1)$ points that does not
depend on the Wasserstein radius and supports a worst-case distribution for every radius. We compute the exact
robust value by checking $\lambda=0$ and at most $Nm$ positive breakpoints, and we construct a worst-case
distribution that moves each observation to a chain corner and splits at most one observation between two
adjacent corners. We apply the results to a convex remediation-planning problem with Monod biodegradation.
Numerical experiments show that the method requires far fewer PDE solves than enumerating all corners. Finally, a
counterexample shows that supermodularity can fail for the fractional $p$-Laplacian when $p>2$.
\end{abstract}

\medskip

\begin{keywords}
distributionally robust optimization, Wasserstein distance, supermodularity, maximum principle, submodular
function minimization, PDE-constrained optimization, fractional Laplacian
\end{keywords}

\begin{MSCcodes}
90C17, 90C27, 90C15, 49Q22, 35J61, 35B50, 35R11
\end{MSCcodes}

\section{Introduction}\label{sec:intro}

Consider an uncertain semilinear elliptic problem
\begin{equation}\label{eq:model}
  -\nabla\!\cdot(A(x)\nabla u)+a(x)u+g(u)=\sum_{k=1}^m\xi_k f_k(x)\quad\text{in }\Omega,\qquad u=0\ \text{on }\partial\Omega.
\end{equation}
Here $\Omega\subset\R^d$ is a bounded domain. The nonnegative function $f_k$ is the spatial profile of source $k$, and $\xi=(\xi_1,\dots,\xi_m)$ gives the uncertain source intensities. The sources may represent pollutant inputs, heat sources, or other nonnegative loads. We measure damage by
\begin{equation}\label{eq:damage}
  F(\xi)=\int_\Omega\Phi\bigl(x,u_\xi(x)\bigr)\dx ,
\end{equation}
where $\Phi(x,\cdot)$ is nondecreasing and convex. Given observations $\xih_1,\dots,\xih_N$ in the box
$\Xi=[\xil,\xib]$, we define the worst-case expected damage by
\begin{equation}\label{eq:dro}
  J(\rho)=\sup\Bigl\{\,\mathbb{E}_P[F(\xi)]\;:\;P\in\Pcal(\Xi),\ W(P,\widehat P_N)\le\rho\Bigr\}.
\end{equation}
The distance $W$ is the type-1 Wasserstein distance associated with the weighted $\ell_1$ norm $\|z\|_w=\sum_kw_k|z_k|$, with
$w>0$.

Previous work on distributionally robust PDE-constrained optimization has used measure approximation \cite{Kouri2017}, approximation schemes \cite{MilzUlbrich2022}, and Rockafellian relaxation \cite{AntilCarneyDiazRoyset2026}. We instead ask when the structure of the PDE makes a Wasserstein distributionally robust optimization problem exactly solvable. Our approach has three steps. The comparison principle gives an order relation between solutions. A concave reaction makes the damage $F$ supermodular. Supermodularity reduces the worst-case optimization problem to a small number of nested stress sets. Wasserstein duality
\cite{MohajerinEsfahaniKuhn2018,BlanchetMurthy2019,GaoKleywegt2023} reduces~\eqref{eq:dro} to inner maximization
problems of the form $\sup_{\xi\in\Xi}\{F(\xi)-\lambda\|\xi-\xih_n\|_w\}$. These problems can be difficult even when $F$ is convex and nondecreasing; in general, maximizing over the box can be NP-hard (\cref{rem:hardness}).

Formal differentiation of the state equation suggests why the result should hold. Write $v_i=\partial u/\partial\xi_i$ and
$w_{ij}=\partial^2u/(\partial\xi_i\partial\xi_j)$. Then
\begin{equation}\label{eq:sens}
  \bigl[-\nabla\!\cdot(A\nabla\cdot)+a+g'(u)\bigr]v_i=f_i,\qquad
  \bigl[-\nabla\!\cdot(A\nabla\cdot)+a+g'(u)\bigr]w_{ij}=-g''(u)\,v_iv_j .
\end{equation}
If $g'\ge0$ and $g''\le0$, the maximum principle suggests that $v_i\ge0$ and $w_{ij}\ge0$. It follows formally that the mixed second derivatives of $F$ are nonnegative:
$\partial^2F/\partial\xi_i\partial\xi_j=\int[\Phi''(u)v_iv_j+\Phi'(u)w_{ij}]\ge0$. However, this argument requires additional smoothness of $g$, $\Phi$, and the solution map. It also excludes nonsmooth damage functions such as $\Phi(r)=(r-u_{\rm safe})_+^2$ and is less suitable for fractional operators. We therefore use comparison arguments that do not require differentiating the solution map.

\paragraph{Contributions}
\begin{enumerate}[leftmargin=*,label=(C\arabic*)]
\item \emph{Order properties.} For coercive bilinear forms satisfying a positivity-preserving lattice condition, we
  prove that the solution operator is nondecreasing when the reaction is nondecreasing and Lipschitz. If the reaction
  is also concave, the solution operator is convex and has increasing differences (\cref{thm:order}). Consequently,
  the damage $F$ is nondecreasing, convex, and supermodular (\cref{thm:damage,cor:fractional}). The proofs use
  comparison and secant-slope arguments and cover local and fractional elliptic equations, equations with advection,
  weighted graphs, and monotone discretizations.
\item \emph{Necessity of concavity.} When the loads are nonnegative, the background load is zero, and $g(0)=0$,
  concavity of $g$ on $[0,\infty)$ is necessary and sufficient for these properties to hold for every operator and
  every damage function in our class (\cref{thm:characterization}). We prove necessity by passing to a large-domain
  limit, in which the PDE response approaches the inverse scalar function $(a\,\mathrm{id}+g)^{-1}$. We also show that
  a strictly convex reaction reverses increasing differences on every domain and for every admissible operator when
  two source profiles coincide and the damage is linear (\cref{thm:convex-reaction}). For a fixed damage function,
  \cref{cor:tradeoff} gives a necessary condition, in the large-domain limit, that relates the curvature of the
  reaction to that of the damage.
\item \emph{Exact robust optimization.} We show that at most $m+1$ nested stress corners per observation are
  needed for the inner maximization problems (\cref{thm:corners,prop:supermodular-set,thm:chain}). The resulting set
  of at most $N(m+1)$ corners does not depend on the Wasserstein radius and supports a worst-case distribution for
  every radius. The robust value follows from checking $\lambda=0$ and at most $Nm$ positive breakpoints
  (\cref{thm:exact}). The worst-case distribution we construct uses at most $N+1$ points: each observation moves to a
  chain corner, and at most one observation is split between two adjacent corners. The bound of at most $N+1$
  support points is already known for general Wasserstein balls \cite{GaoKleywegt2023}; our contribution is to
  locate these points on the stress chains (\cref{rem:novelty}). The results extend to weighted observations and
  observation-dependent transport costs (\cref{rem:weighted}). A chain can be computed with at most $2m-1$ submodular
  minimizations per observation (\cref{prop:algorithm}), and the chain ranks the sources by their importance to the
  worst-case outcome (\cref{rem:ranking}). In tests with $m=20$ sources, the method used about $960$ times fewer PDE
  solves than full corner enumeration.
\item \emph{Application and limits.} We apply the results to remediation-capacity planning with Monod
  biodegradation. The resulting robust design problem is convex and can be solved by a cutting-plane method using
  worst-case subgradients (\cref{sec:application,sec:numerics}). Synthetic experiments independently check the robust
  values and compare the number of PDE solves with full corner enumeration (\cref{sec:numerics}). The stress chains
  agree across the tested grids, and the optimal design also agrees on the grids for which it was computed. Finally,
  \cref{sec:plaplace} gives a counterexample to supermodularity for the fractional $p$-Laplacian when $p>2$.
\end{enumerate}

\paragraph{Related work} Three areas of research are most relevant to this paper.

\emph{Wasserstein distributionally robust optimization and PDE constraints.} Wasserstein duality and
worst-case distributions are studied in \cite{MohajerinEsfahaniKuhn2018,BlanchetMurthy2019,GaoKleywegt2023}.
In particular, for a Wasserstein ball around an empirical distribution with $N$ atoms, a worst-case distribution
supported on at most $N+1$ points, obtained by moving each observation and splitting at most one of them, exists
under general conditions \cite{GaoKleywegt2023}. Our \cref{thm:exact} recovers this cardinality; its contribution is
to locate the support points on stress chains generated by the PDE, in a set that does not depend on the radius,
and to express the dual function using finitely many breakpoints.
Distributionally robust problems with PDE constraints are studied in
\cite{Kouri2017,MilzUlbrich2022,AntilCarneyDiazRoyset2026}. These papers provide general methods for handling
uncertain probability distributions. We instead identify PDE conditions that reduce the robust problem exactly to a finite collection of sets.

\emph{Supermodularity in optimization.} Supermodularity and monotone comparative statics are standard tools in
operations research \cite{Topkis1998}. For distributions with prescribed marginals, comonotone couplings
maximize expectations of supermodular functions \cite{Tchen1980}, a fact used to study the price of correlations
\cite{AgrawalDingSaberiYe2012}. Long, Qi, and Zhang \cite{LongQiZhang2024} use supermodularity in two-stage
distributionally robust optimization. Leng et al.\ \cite{LengXieRanZhang2025} use it to reduce disruption scenarios
in supply chain network design. Robust and distributionally robust submodular maximization are studied in
\cite{KrauseEtAl2008,StaibWilderJegelka2019}. In these studies, supermodularity arises from the stochastic,
recourse, or discrete optimization model. In our setting, it follows from the elliptic PDE and the reaction term,
and we give a necessary and sufficient condition for it. We then use it, under weighted $\ell_1$ Wasserstein
uncertainty sets, to compute the exact worst-case value and to construct a worst-case distribution.

\emph{Parametric submodular minimization.} For a submodular function $f$ and a nondecreasing modular function $c$,
the smallest and largest minimizers of $f+\lambda c$ are nested as $\lambda$ changes \cite{Topkis1978}. This result
is related to principal partitions and strong map sequences \cite{Fujishige2005,IwataMurotaShigeno1997}, as well
as to parametric flow and submodular algorithms \cite{GalloGrigoriadisTarjan1989,FleischerIwata2003}. Our
\cref{thm:chain} uses this known property. Our contribution is to derive the parametric submodular problem from the PDE and to show that its nested sets determine the robust value and a worst-case distribution. Submodular
minimization can be performed in strongly polynomial time \cite{Schrijver2000,IwataFleischerFujishige2001,Orlin2009}.
The comparison arguments in \cref{sec:order} are based on standard results for elliptic operators
\cite{GilbargTrudinger2001} and Dirichlet forms \cite{FukushimaOshimaTakeda2011}. We combine these tools to obtain an exact solution method for the robust problem.

\section{Mathematical framework}\label{sec:framework}

\subsection{Standing assumptions}
Let $(X,\mu)$ be a finite measure space, and order $L^2:=L^2(X,\mu)$ by almost-everywhere comparison. For a function $v$, we
write $v^+=\max\{v,0\}$ and $v^-=\max\{-v,0\}$, so that $v=v^+-v^-$.

\begin{assumption}[Operator]\label[assumption]{ass:H}
\leavevmode
\begin{enumerate}[label=(H\arabic*),leftmargin=*]
\item\label{H1} $V$ is a real Hilbert space, continuously and densely embedded in $L^2$, with
  $\|v\|_{L^2}\le C_V\|v\|_V$.
\item\label{H2} $\E:V\times V\to\R$ is bilinear, bounded ($|\E(v,\varphi)|\le M_\E\|v\|_V\|\varphi\|_V$) and
  coercive: $\E(v,v)\ge\alpha\|v\|_V^2$ for some $\alpha>0$.
\item\label{H3} (positivity-preserving lattice property) For every $v\in V$, we have $v^+\in V$ and $\E(v^+,v^-)\le0$.
\end{enumerate}
\end{assumption}

Condition~\ref{H3} is related to the first Beurling--Deny criterion \cite{FukushimaOshimaTakeda2011}, but it
does not require the stronger Markov (unit-contraction) property. Define $\Vp=\{\varphi\in V:\varphi\ge0\}$. We regard $h\in L^2$ as an element of $V^*$ through the functional $\varphi\mapsto\int_Xh\varphi\dmu$. We order $V^*$ by the positive
cone $\Vsp=\{\ell\in V^*:\ip{\ell}{\varphi}\ge0\ \forall\varphi\in\Vp\}$; that is, $\ell\le\ell'$ means $\ell'-\ell\in\Vsp$.

\begin{assumption}[Reaction]\label[assumption]{ass:G}
$g:\R\to\R$ is nondecreasing and Lipschitz with constant $L_g$.
\end{assumption}

\begin{assumption}[Damage density]\label[assumption]{ass:Phi}
The function $\Phi:X\times\R\to\R$ satisfies the following: $\Phi(\cdot,r)$ is measurable for each $r$;
$\Phi(x,\cdot)$ is nondecreasing and convex for a.e.\ $x$; and $|\Phi(x,r)|\le a_0(x)+b_0r^2$ with $a_0\in L^1(X)$
and $b_0\ge0$.
\end{assumption}

For a given load $\ell\in V^*$, we seek $u\in V$ satisfying
\begin{equation}\label{eq:state}
  \E(u,\varphi)+\int_Xg(u)\varphi\dmu=\ip{\ell}{\varphi}\qquad\forall\varphi\in V.
\end{equation}
The bound $|g(r)|\le|g(0)|+L_g|r|$ and the finiteness of $\mu(X)$ imply that $g(u)\in L^2$ whenever $u\in L^2$.
Thus~\eqref{eq:state} is well defined. We represent uncertain loads by
\begin{equation}\label{eq:loads}
  \ell(\xi)=\ell_0+\sum_{k=1}^m\xi_kf_k,\qquad \ell_0\in V^*,\quad f_1,\dots,f_m\in\Vsp,\quad \xi\in\R^m.
\end{equation}
The fixed load $\ell_0$ need not be nonnegative. In \cref{sec:application}, it includes the effect of remediation.

\subsection{Examples}
The following four examples satisfy \ref{H1}--\ref{H3}, as verified in \cref{app:examples}.

\begin{example}[Diffusion--advection--reaction]\label[example]{ex:local}
Let $\Omega\subset\R^d$ be a bounded Lipschitz domain, let $X=\Omega$ with Lebesgue measure, let
$V=H^1_0(\Omega)$, and set
\[
  \E(u,\varphi)=\int_\Omega\bigl(A\nabla u\cdot\nabla\varphi+(\beta\cdot\nabla u)\varphi+au\varphi\bigr)\dx,
\]
where $A\in L^\infty(\Omega;\R^{d\times d})$ is uniformly elliptic, $a\in L^\infty$, $\beta\in W^{1,\infty}(\Omega;\R^d)$,
and $a-\tfrac12\operatorname{div}\beta\ge0$ (for example, $a\ge0$ and $\operatorname{div}\beta=0$).
\end{example}

\begin{example}[Fractional Laplacian]\label[example]{ex:fractional}
Let $\Omega\subset\R^d$ be a bounded open set, let $s\in(0,1)$ and $X=\Omega$, let
$V=X^s_0(\Omega):=\{u\in H^s(\R^d):u=0\text{ a.e.\ in }\R^d\setminus\Omega\}$ with the Gagliardo seminorm, and set
\[
  \E(u,\varphi)=\frac{C_{d,s}}2\iint_{\R^{2d}}\frac{(u(x)-u(y))(\varphi(x)-\varphi(y))}{|x-y|^{d+2s}}\dx\,\mathrm{d}y
  +\int_\Omega au\varphi\dx,\qquad a\in L^\infty(\Omega),\ a\ge0 .
\]
This is the weak formulation of $(-\Delta)^su+au$ with zero exterior Dirichlet data \cite{ServadeiValdinoci2013}.
\end{example}

\begin{example}[Weighted graphs and networks]\label[example]{ex:graph}
Let $X$ be a finite set of interior nodes with positive weights $\mu_x$, let $V=\R^X$, and set
$\E(u,\varphi)=\frac12\sum_{x,y}\omega_{xy}(u_x-u_y)(\varphi_x-\varphi_y)+\sum_x\kappa_xu_x\varphi_x$, where
$\omega_{xy}=\omega_{yx}\ge0$, $\kappa_x\ge0$, and $\kappa\neq0$ on every connected component (absorption, or a
connection to boundary nodes).
\end{example}

\begin{example}[Monotone discretizations]\label[example]{ex:discrete}
Let $X=\{1,\dots,n\}$ with nodal weights $\mu_i>0$, let $V=\R^n$, and set $\E(u,\varphi)=\varphi^\top Ku$, where
$K$ is a Z-matrix ($K_{ij}\le0$ for $i\ne j$) and $K+K^\top$ is positive definite. The reaction is applied at the
nodes with lumped weights, $\int g(u)\varphi\dmu=\sum_i\mu_ig(u_i)\varphi_i$, and a load $f\in\R^n$ acts through
$\varphi\mapsto\sum_i\mu_if_i\varphi_i$. If $K^{\rm FD}$ is a finite-difference matrix for the strong form of the
operator and $M=\operatorname{diag}(\mu_1,\dots,\mu_n)$, we take $K=MK^{\rm FD}$; then~\eqref{eq:state} reads
$K^{\rm FD}u+g(u)=f$ at the nodes. For uniform weights, $K$ is a Z-matrix with positive definite symmetric part
exactly when $K^{\rm FD}$ is. This example applies to discretizations whose assembled matrix is a
Z-matrix with positive definite symmetric part. Examples include the five-point and upwind finite-difference schemes
used in \cref{sec:numerics}; $P_1$ finite elements with mass lumping for $-\Delta$ with a nonnegative reaction on
non-obtuse meshes \cite{CiarletRaviart1973}; the monotone finite element scheme of \cite{XuZikatanov1999} for
convection--diffusion, under the mesh conditions given there; and monotone schemes for the fractional Laplacian
\cite{HuangOberman2014}. A standard $P_1$ advection discretization need not produce a Z-matrix. Using a full rather than lumped mass matrix may also destroy this property.
\end{example}

\subsection{Three lemmas}

\begin{lemma}[Well-posedness]\label[lemma]{lem:wellposed}
Assume \ref{H1}--\ref{H2} and \cref{ass:G}. For every $\ell\in V^*$, problem~\eqref{eq:state} has a
unique solution $u=\mathcal S\ell$, and $\|\mathcal S\ell-\mathcal S\ell'\|_V\le\alpha^{-1}\|\ell-\ell'\|_{V^*}$ for all $\ell,\ell'\in V^*$.
\end{lemma}
\begin{proof}
Define $\mathcal A:V\to V^*$ by $\ip{\mathcal Av}{\varphi}=\E(v,\varphi)+\int g(v)\varphi\dmu$. Since $g$ is
nondecreasing, $\ip{\mathcal Av-\mathcal Aw}{v-w}\ge\E(v-w,v-w)\ge\alpha\|v-w\|_V^2$, and
$\|\mathcal Av-\mathcal Aw\|_{V^*}\le(M_\E+L_gC_V^2)\|v-w\|_V$. Thus $\mathcal A$ is strongly monotone and Lipschitz continuous. Its bijectivity follows from the main theorem on strongly monotone operators
\cite[Thm.~25.B]{Zeidler1990}. For completeness, we also give a direct argument. Let $R:V^*\to V$ be the Riesz map and $L_{\mathcal A}=M_\E+L_gC_V^2$. Then
$\|v-w-tR(\mathcal Av-\mathcal Aw)\|_V^2\le(1-2t\alpha+t^2L_{\mathcal A}^2)\|v-w\|_V^2$, so
$v\mapsto v-tR(\mathcal Av-\ell)$ is a contraction on $V$ for $0<t<2\alpha/L_{\mathcal A}^2$, and its unique fixed
point solves~\eqref{eq:state}. For $u_1=\mathcal S\ell$ and $u_2=\mathcal S\ell'$, testing the difference of the two equations with $u_1-u_2$ gives $\alpha\|u_1-u_2\|_V^2\le\ip{\ell-\ell'}{u_1-u_2}$, which yields the
estimate.
\end{proof}

\begin{lemma}[Comparison]\label[lemma]{lem:comparison}
Assume \ref{H1}--\ref{H3}. Let $c\in L^\infty(X)$, $c\ge0$, and let $v\in V$ satisfy
$\E(v,\varphi)+\int cv\varphi\dmu\ge0$ for all $\varphi\in\Vp$. Then $v\ge0$ a.e.
\end{lemma}
\begin{proof}
By \ref{H3}, $v^-=v^+-v\in\Vp$. Testing with $v^-$ and using $vv^-=-(v^-)^2$, we obtain
\[
0\le\E(v,v^-)+\int cvv^-\dmu=\E(v^+,v^-)-\E(v^-,v^-)-\int c(v^-)^2\dmu\le-\alpha\|v^-\|_V^2 .
\]
Thus $v^-=0$.
\end{proof}

For $p,q\in\R$, we define the secant slope
\[
  s_g(p,q)=\frac{g(q)-g(p)}{q-p}\ (p\ne q),\qquad
  s_g(p,p)=\begin{cases}g'_+(p)&\text{if $g$ is concave or convex},\\ 0&\text{otherwise.}\end{cases}
\]

\begin{lemma}[Secant slopes]\label[lemma]{lem:secant}
Let $g$ satisfy \cref{ass:G}. Then $s_g$ is Borel measurable and symmetric, it takes values in $[0,L_g]$, and
$g(q)-g(p)=s_g(p,q)(q-p)$. Moreover:
\begin{enumerate}[label=(\roman*)]
\item if $g$ is concave, then $s_g(p,q)\ge s_g(p',q')$ whenever $p\le p'$ and $q\le q'$;
\item if $g$ is convex, then $s_g(p,q)\le s_g(p',q')$ whenever $p\le p'$ and $q\le q'$; if $g$ is strictly
  convex and $(p,q)\ne(p',q')$, the inequality is strict.
\end{enumerate}
\end{lemma}
\begin{proof}
Measurability and the bounds follow from the definition. Suppose first that $g$ is concave. The three-chord
inequality gives $s_g(x,y)\ge s_g(x,z)\ge s_g(y,z)$ for $x<y<z$. In addition, $s_g(p,q)\le g'_+(p)$ for $q>p$ and
$s_g(p,q)\ge g'_-(p)\ge g'_+(p)$ for $q<p$. Thus $q\mapsto s_g(p,q)$ is nonincreasing; by symmetry, so is
$p\mapsto s_g(p,q)$. Therefore $s_g(p,q)\ge s_g(p',q)\ge s_g(p',q')$, which proves (i). Part (ii) follows by applying
(i) to $-g$. If $g$ is strictly convex, every inequality used above is strict.
\end{proof}

\section{Order structure of the solution operator}\label{sec:order}

\begin{theorem}[Order structure]\label[theorem]{thm:order}
Assume \ref{H1}--\ref{H3} and \cref{ass:G}, and let $\mathcal S:V^*\to V$ be the solution operator
of~\eqref{eq:state}.
\begin{enumerate}[label=(\alph*)]
\item\label{it:mono} (Monotonicity) If $\ell\le\ell'$, then $\mathcal S\ell\le\mathcal S\ell'$.
\item\label{it:convex} (Convexity) If $g$ is concave, then for all $\ell,\ell'\in V^*$ and $\theta\in[0,1]$, \[\mathcal S\bigl(\theta\ell+(1-\theta)\ell'\bigr)\le\theta\,\mathcal S\ell+(1-\theta)\,\mathcal S\ell'.\]
\item\label{it:incdiff} (Increasing differences) If $g$ is concave, then for all $\ell\in V^*$ and $p,q\in\Vsp$,
  \[\mathcal S(\ell+p+q)-\mathcal S(\ell+q)\ \ge\ \mathcal S(\ell+p)-\mathcal S(\ell)\ \ge\ 0.\]
\item\label{it:reverse} (Reversal) If $g$ is convex, then for all $\ell,\ell'\in V^*$, $\theta\in[0,1]$, and $p,q\in\Vsp$, \[\mathcal S\bigl(\theta\ell+(1-\theta)\ell'\bigr)\ge\theta\,\mathcal S\ell+(1-\theta)\,\mathcal S\ell',\qquad
  0\le\mathcal S(\ell+p+q)-\mathcal S(\ell+q)\le\mathcal S(\ell+p)-\mathcal S(\ell).\]
\end{enumerate}
\end{theorem}

\begin{proof}
\ref{it:mono} Let $u_1=\mathcal S\ell$ and $u_2=\mathcal S\ell'$. We subtract the equations for $u_1$ and $u_2$. With $d=u_2-u_1$ and
$c=s_g(u_1,u_2)\in[0,L_g]$, we have $g(u_2)-g(u_1)=cd$, so the difference satisfies
$\E(d,\varphi)+\int cd\varphi=\ip{\ell'-\ell}{\varphi}\ge0$ for $\varphi\in\Vp$. \Cref{lem:comparison} gives $d\ge0$.

\ref{it:convex} Let $u_1=\mathcal S\ell$, $u_2=\mathcal S\ell'$, $w=\theta u_1+(1-\theta)u_2$, $\zeta=\theta\ell+(1-\theta)\ell'$, and
$u=\mathcal S\zeta$. By linearity, $\E(w,\varphi)+\int[\theta g(u_1)+(1-\theta)g(u_2)]\varphi=\ip{\zeta}{\varphi}$.
Since $g$ is concave, $\theta g(u_1)+(1-\theta)g(u_2)\le g(w)$ pointwise, so $w$ is a supersolution for the load
$\zeta$: for $\varphi\in\Vp$,
\[
  \E(w,\varphi)+\int g(w)\varphi\dmu\ \ge\ \ip{\zeta}{\varphi}=\E(u,\varphi)+\int g(u)\varphi\dmu .
\]
With $d=w-u$ and $c=s_g(u,w)$, this reads $\E(d,\varphi)+\int cd\varphi\ge0$. Applying \cref{lem:comparison}
to $d$ gives $w\ge u$.

\ref{it:incdiff} Let $u_{00}=\mathcal S\ell$, $u_{10}=\mathcal S(\ell+p)$, $u_{01}=\mathcal S(\ell+q)$, and
$u_{11}=\mathcal S(\ell+p+q)$, and set $d=u_{11}-u_{01}$ and $e=u_{10}-u_{00}$. By \ref{it:mono}, $e\ge0$,
$u_{01}\ge u_{00}$, and $u_{11}\ge u_{10}$. With $c_1=s_g(u_{01},u_{11})$ and $c_2=s_g(u_{00},u_{10})$, we have
\[
  \E(d,\varphi)+\int c_1d\varphi=\ip{p}{\varphi},\qquad \E(e,\varphi)+\int c_2e\varphi=\ip{p}{\varphi}.
\]
Subtracting the two equations gives
\begin{equation}\label{eq:de}
  \E(d-e,\varphi)+\int c_1(d-e)\varphi\dmu=\int(c_2-c_1)e\varphi\dmu .
\end{equation}
Both arguments defining $c_1$ are at least the corresponding arguments defining $c_2$, so \cref{lem:secant}(i) gives
$c_1\le c_2$. The right-hand side of~\eqref{eq:de} is therefore nonnegative for $\varphi\in\Vp$, and
\cref{lem:comparison} yields $d\ge e$.

\ref{it:reverse} If $g$ is convex, the comparison in part~\ref{it:convex} is reversed:
$\theta g(u_1)+(1-\theta)g(u_2)\ge g(w)$ makes $w$ a subsolution, and \cref{lem:comparison} applied to $u-w$ gives
$u\ge w$. In~\eqref{eq:de}, we now have $c_1\ge c_2$ by \cref{lem:secant}(ii). Writing~\eqref{eq:de} for $e-d$
gives $\E(e-d,\varphi)+\int c_1(e-d)\varphi=\int(c_1-c_2)e\varphi\ge0$, and \cref{lem:comparison} yields $e\ge d$.
\end{proof}

\begin{corollary}[Parametrized solution map]\label[corollary]{cor:param}
In the setting of \cref{thm:order}, let $u_\xi=\mathcal S\ell(\xi)$ with $\ell(\xi)$ as in~\eqref{eq:loads}.
Then $\xi\mapsto u_\xi$ is nondecreasing on $\R^m$. If $g$ is concave, this map is convex and supermodular:
\[
  u_{\xi\vee\eta}+u_{\xi\wedge\eta}\ \ge\ u_\xi+u_\eta\quad\text{a.e., for all }\xi,\eta\in\R^m .
\]
If $g$ is convex, the map is concave and submodular.
\end{corollary}
\begin{proof}
Monotonicity and convexity follow from \cref{thm:order}. To prove supermodularity, we apply
\cref{thm:order}\ref{it:incdiff} with
\[
  \ell=\ell(\xi\wedge\eta),\qquad p=\sum_k(\xi-\xi\wedge\eta)_kf_k,\qquad q=\sum_k(\eta-\xi\wedge\eta)_kf_k,
\]
which lie in $\Vsp$. Since $\xi\vee\eta=\xi+\eta-\xi\wedge\eta$, we have
$\ell+p+q=\ell(\xi\vee\eta)$.
\end{proof}

\begin{lemma}[Two-point inequality]\label[lemma]{lem:twopoint}
Let $\phi:\R\to\R$ be nondecreasing and convex, and let $a,b,c,d\in\R$ satisfy $a\le\min\{b,c\}$ and
$a+d\ge b+c$. Then $\phi(a)+\phi(d)\ge\phi(b)+\phi(c)$.
\end{lemma}
\begin{proof}
Let $d'=b+c-a\le d$. Then $\phi(d)\ge\phi(d')$ and $a\le b,c\le d'$. If $a=d'$, all four numbers coincide.
Otherwise, $\lambda=(d'-b)/(d'-a)\in[0,1]$ satisfies $b=\lambda a+(1-\lambda)d'$ and $c=(1-\lambda)a+\lambda d'$.
Adding the two convexity inequalities gives $\phi(b)+\phi(c)\le\phi(a)+\phi(d')\le\phi(a)+\phi(d)$.
\end{proof}

\begin{theorem}[Damage functional]\label[theorem]{thm:damage}
Assume \ref{H1}--\ref{H3}, \cref{ass:G} with $g$ concave, \cref{ass:Phi}, and~\eqref{eq:loads}. Then
\[
  F(\xi)=\int_X\Phi\bigl(x,u_\xi(x)\bigr)\dmu(x)
\]
is finite, nondecreasing, jointly convex, and supermodular on $\R^m$. In particular, it is continuous and
Lipschitz on bounded sets. If $\Phi(x,\cdot)$ is strictly increasing for a.e.\ $x$ and the $f_k$ are nonzero, then
$F$ is strictly increasing in each coordinate.
\end{theorem}
\begin{proof}
The function $x\mapsto\Phi(x,u_\xi(x))$ is measurable, because $\Phi$ is a Carath\'eodory function: it is
measurable in $x$ and, being convex and finite, continuous in $r$ for a.e.\ $x$. The growth bound in \cref{ass:Phi}
and $u_\xi\in L^2$ imply that $F(\xi)$ is finite. Monotonicity follows from
\cref{cor:param} and the monotonicity of $\Phi(x,\cdot)$. For convexity, let $\zeta=\theta\xi+(1-\theta)\eta$; then
\[
  \Phi(x,u_\zeta)\le\Phi\bigl(x,\theta u_\xi+(1-\theta)u_\eta\bigr)\le\theta\Phi(x,u_\xi)+(1-\theta)\Phi(x,u_\eta),
\]
by \cref{cor:param}, the monotonicity of $\Phi(x,\cdot)$, and then its convexity. To prove supermodularity, we
apply \cref{lem:twopoint} pointwise for almost every $x$ with $a=u_{\xi\wedge\eta}(x)$, $b=u_\xi(x)$, $c=u_\eta(x)$, and
$d=u_{\xi\vee\eta}(x)$. Here $a\le\min\{b,c\}$ by monotonicity and $a+d\ge b+c$ by \cref{cor:param}. Integrating
over $X$ gives the claim. Since $F$ is a finite convex function on $\R^m$, it is continuous and locally Lipschitz.
For strict monotonicity, let $\xi\le\eta$ with $\xi_k<\eta_k$. Then $u_\eta-u_\xi\ge0$ and $u_\eta-u_\xi\ne0$;
otherwise $\ell(\eta)=\ell(\xi)$ by~\eqref{eq:state}, which forces $f_k=0$. Hence $\Phi(x,u_\eta)>\Phi(x,u_\xi)$ on a
set of positive measure.
\end{proof}

\begin{corollary}[Fractional elliptic loads]\label[corollary]{cor:fractional}
Let $\Omega\subset\R^d$ be bounded and open, and let $s\in(0,1)$. Let $a\in L^\infty(\Omega)$ with $a\ge0$, let $g$
be nondecreasing, concave, and Lipschitz, let $f_k\in L^2(\Omega)$ with $f_k\ge0$, and let $\Phi$ be as in
\cref{ass:Phi}. Let $u_\xi$ be the weak solution of
\[
  (-\Delta)^su+a(x)u+g(u)=\sum_{k=1}^m\xi_kf_k\ \text{ in }\Omega,\qquad u=0\ \text{ in }\R^d\setminus\Omega .
\]
Then $F(\xi)=\int_\Omega\Phi(x,u_\xi)\dx$ is nondecreasing, jointly convex, and supermodular in $\xi\in\R^m$.
\end{corollary}
\begin{proof}
The operator in \cref{ex:fractional} satisfies \ref{H1}--\ref{H3} (\cref{app:examples}), so the result
follows from \cref{thm:damage}.
\end{proof}

\section{Necessity of concavity}\label{sec:sharp}

\subsection{Convex reactions reverse the order properties}

\begin{theorem}[Strict failure for convex reactions]\label[theorem]{thm:convex-reaction}
Assume \ref{H1}--\ref{H3} and \cref{ass:G} with $g$ strictly convex. Let $f_1=f_2=:f\in\Vsp\setminus\{0\}$, let
$\ell_0\in V^*$ and $f_3,\dots,f_m\in\Vsp$ be arbitrary, and take $\Phi(x,r)=r$. Then, for every $\xi^0\in\R^m$
and $t>0$,
\begin{equation}\label{eq:strict}
  F(\xi^0+te_1+te_2)+F(\xi^0)\ <\ F(\xi^0+te_1)+F(\xi^0+te_2).
\end{equation}
Thus, when the reaction is strictly convex, $F$ fails to have increasing differences in $(\xi_1,\xi_2)$ at every
point. This holds for every domain and every operator satisfying \cref{ass:H}.
\end{theorem}
\begin{proof}
Let $u_0=u_{\xi^0}$, $u_1=u_{\xi^0+te_1}=u_{\xi^0+te_2}$ (the two loads coincide), and
$u_2=u_{\xi^0+te_1+te_2}$. Set $d=u_2-u_1$ and $e=u_1-u_0$. Both solve linear equations with load $tf$ and
coefficients $c_1=s_g(u_1,u_2)$ and $c_2=s_g(u_0,u_1)$, so~\eqref{eq:de} holds. By
\cref{thm:order}\ref{it:mono}, $e\ge0$; moreover, $e\ne0$, since otherwise $tf=0$. On $\{e>0\}$, we have
$u_1>u_0$ and $u_2\ge u_1$, so \cref{lem:secant}(ii) gives $c_1>c_2$ there. Hence $\chi:=(c_2-c_1)e\le0$, with
$\chi<0$ on a set of positive measure. Applying \cref{lem:comparison} to $e-d$ gives $d-e\le0$. Moreover,
$d-e\ne0$; otherwise $\int\chi\varphi=0$ for all $\varphi$ in the dense subspace $V$, which forces $\chi=0$.
Therefore $F(\xi^0+te_1+te_2)+F(\xi^0)-2F(\xi^0+te_1)=\int(d-e)\dmu<0$.
\end{proof}

\begin{example}\label[example]{ex:softplus}
The softplus function $g(r)=\kappa(\log(1+e^r)-\log2)$ is smooth, increasing, strictly convex, and
$\kappa$-Lipschitz, with $g(0)=0$. Therefore, \cref{thm:convex-reaction} applies on any domain and to every
operator in \cref{ex:local,ex:fractional,ex:graph,ex:discrete}.
\end{example}

\subsection{A characterization}

We now show that concavity is necessary when the result is required to hold for every admissible problem. Let $Q_R=(-R,R)^d$ and, for fixed $a>0$, consider the Dirichlet problem
\begin{equation}\label{eq:cube}
  -\Delta u+au+g(u)=\xi_1+\xi_2\ \text{ in }Q_R,\qquad u=0\ \text{ on }\partial Q_R,
\end{equation}
which is an instance of \cref{ex:local} with $f_1=f_2=1$. Set $h(r)=ar+g(r)$. When $g(0)=0$, $h$ is an increasing
bijection of $[0,\infty)$ onto itself, and we write $H=h^{-1}$ for its inverse.

\begin{lemma}[Large-cube asymptotics]\label[lemma]{lem:cube}
Let $g$ satisfy \cref{ass:G} with $g(0)=0$, let $t\ge0$, and let $u_R$ solve~\eqref{eq:cube} with
$\xi_1+\xi_2=t$. Then
\[
  0\le u_R\le H(t)\quad\text{and}\quad
  0\le H(t)-\frac1{|Q_R|}\int_{Q_R}u_R\dx\le\frac{d\,H(t)}{\sqrt a\,R}.
\]
\end{lemma}
\begin{proof}
The constant $U=H(t)\ge0$ solves the equation in~\eqref{eq:cube} inside $Q_R$; since its boundary values $U\ge0$ lie above the zero Dirichlet data, it is a supersolution of the Dirichlet problem. We use the weak maximum principle for $-\Delta+c$ with
$c\ge0$ bounded and with boundary inequalities in the $H^1$ sense \cite[Thm.~8.1]{GilbargTrudinger2001}. Since the
load is nonnegative and $\mathcal S0=0$, \cref{thm:order}\ref{it:mono} gives $u_R\ge0$. The function $U-u_R$ is
nonnegative on $\partial Q_R$ and satisfies $-\Delta(U-u_R)+(a+s_g(u_R,U))(U-u_R)=0$. By comparison,
$0\le u_R\le U$. Next, let $\vartheta(x)=\sum_{i=1}^d\cosh(\sqrt a\,x_i)/\cosh(\sqrt a\,R)$, which satisfies
$-\Delta\vartheta+a\vartheta=0$ and $\vartheta\ge1$ on $\partial Q_R$. Put $z=u_R-U(1-\vartheta)$; then $z\ge0$ on
$\partial Q_R$, and since $u_R\le U$ and $g$ is nondecreasing,
\[
  -\Delta z+az=t-g(u_R)-aU\ \ge\ t-g(U)-aU=0
\]
weakly. Therefore $z\ge0$; that is, $U-u_R\le U\vartheta$. Finally,
\[
  |Q_R|^{-1}\int_{Q_R}\vartheta\dx=\frac{d\tanh(\sqrt aR)}{\sqrt aR}\le\frac{d}{\sqrt aR}.
\]
\end{proof}

\begin{theorem}[Characterization]\label[theorem]{thm:characterization}
Let $g$ satisfy \cref{ass:G} with $g(0)=0$. The following are equivalent:
\begin{enumerate}[label=(\roman*)]
\item\label{ch1} $g$ is concave on $[0,\infty)$.
\item\label{ch2} For every instance of \cref{ass:H}, every $f_1,\dots,f_m\in\Vsp$ with $\ell_0=0$, and every $\Phi$
  satisfying \cref{ass:Phi}, the damage $F$ is nondecreasing, convex, and supermodular on $\R^m_+$.
\item\label{ch3} For some $a>0$ and some sequence $R_j\to\infty$, the normalized total mass
  $F_j(\xi)=|Q_{R_j}|^{-1}\int_{Q_{R_j}}u_\xi\dx$ of~\eqref{eq:cube} is supermodular on $\R^2_+$.
\end{enumerate}
\end{theorem}
\begin{proof}
\ref{ch1}$\Rightarrow$\ref{ch2}. Extend $g$ to the negative half-line by $\tilde g=g$ on $[0,\infty)$ and
$\tilde g(r)=g'_+(0)\,r$ for $r<0$. The extension $\tilde g$ is nondecreasing, concave, and Lipschitz on $\R$. For
$\xi\in\R^m_+$, \cref{thm:order}\ref{it:mono} shows that the $\tilde g$-solution is at least $\tilde{\mathcal S}0=0$.
It therefore also solves the problem with $g$, and by uniqueness the two solution maps agree on $\R^m_+$. Applying
\cref{thm:damage} to $\tilde g$ and restricting to the sublattice $\R^m_+$ gives the result.

\ref{ch2}$\Rightarrow$\ref{ch3}. Take \cref{ex:local} on $Q_R$ with $A=I$, $\beta=0$, constant $a>0$, $f_1=f_2=1$,
and $\Phi(x,r)=r/|Q_R|$.

\ref{ch3}$\Rightarrow$\ref{ch1}. Since the two source profiles are identical, $F_j$ depends only on $\xi_1+\xi_2$. For $T\ge0$ and $\delta>0$,
supermodularity at $\xi=(T+\delta,0)$ and $\eta=(T,\delta)$ gives
$F_j(T+\delta,\delta)+F_j(T,0)\ge F_j(T+\delta,0)+F_j(T,\delta)$. Letting $j\to\infty$ and using \cref{lem:cube}, we
obtain the midpoint-convexity inequality $H(T+2\delta)+H(T)\ge2H(T+\delta)$. Since $H$ is continuous, it is convex
on $[0,\infty)$. Therefore its increasing inverse $h$ is concave, and so is $g=h-a\,\mathrm{id}$ on $[0,\infty)$.
\end{proof}

\begin{corollary}[Curvature trade-off for a fixed damage]\label[corollary]{cor:tradeoff}
Let $g$ be as in \cref{thm:characterization}, let $\Phi:[0,\infty)\to\R$ be nondecreasing and convex, and suppose that
$\xi\mapsto|Q_{R_j}|^{-1}\int_{Q_{R_j}}\Phi(u_\xi)\dx$ is supermodular on $\R^2_+$ along a sequence $R_j\to\infty$.
Then $\Phi\circ H$ is convex on $[0,\infty)$. If $g,\Phi\in C^2$ and $\Phi'>0$, this is equivalent to
\begin{equation}\label{eq:tradeoff}
  \frac{\Phi''(r)}{\Phi'(r)}\ \ge\ \frac{g''(r)}{a+g'(r)}\qquad(r\ge0).
\end{equation}
\end{corollary}
\begin{proof}
On $[0,U]$, the convex nondecreasing function $\Phi$ is Lipschitz with constant $\Phi'_-(U)$. Hence
$|\Phi(u_R)-\Phi(U)|\le\Phi'_-(U)(U-u_R)$, and \cref{lem:cube} gives $|Q_R|^{-1}\int\Phi(u_R)\to\Phi(H(t))$. Passing to the large-domain limit as in the proof of \cref{thm:characterization} shows that $\Phi\circ H$ is convex.
Differentiating gives $(\Phi\circ H)''=(\Phi''h'-\Phi'h'')/h'^3$ at $r=H(t)$, with $h'=a+g'>0$, which
yields~\eqref{eq:tradeoff}.
\end{proof}

\section{Exact robust optimization using stress corners}\label{sec:dro}

In this section, let $F:\R^m\to\R$ be continuous, nondecreasing, convex, and supermodular. \Cref{thm:damage} shows
that these assumptions hold for the damage functions considered above when the reaction is concave. Let
$\Xi=[\xil,\xib]\subset\R^m$ be a box, let $\xih_1,\dots,\xih_N\in\Xi$ be observations, and let
$\widehat P_N=\frac1N\sum_n\delta_{\xih_n}$ and $\|z\|_w=\sum_kw_k|z_k|$ with $w_k>0$. The type-1 Wasserstein
distance is $W(P,Q)=\inf_{\gamma\in\Gamma(P,Q)}\int\|\xi-\zeta\|_w\,\mathrm d\gamma$, and $J(\rho)$ is defined
by~\eqref{eq:dro} with this $W$.

\begin{definition}[Stress corners]\label[definition]{def:corners}
For $n\le N$ and $S\subseteq[m]$, let $\xih_n^S$ be the point with coordinates $\xib_k$ for $k\in S$ and
$\xih_{nk}$ for $k\notin S$. Define
\begin{gather*}
  c_n(S)=\|\xih_n^S-\xih_n\|_w=\sum_{k\in S}w_k(\xib_k-\xih_{nk}),\\
  \Psi_n(S;\lambda)=F(\xih_n^S)-\lambda c_n(S),\qquad
  \varphi_n(\lambda)=\max_{S\subseteq[m]}\Psi_n(S;\lambda).
\end{gather*}
Let $K_n=\{k:\xih_{nk}=\xib_k\}$ be the set of coordinates that are already at their upper bound.
\end{definition}

\begin{theorem}[Upper stress corners]\label[theorem]{thm:corners}
For every $n$ and $\lambda\ge0$,
\[
  \sup_{\xi\in\Xi}\bigl\{F(\xi)-\lambda\|\xi-\xih_n\|_w\bigr\}=\varphi_n(\lambda),
\]
and the supremum is attained at an upper stress corner $\xih_n^S$.
\end{theorem}
\begin{proof}
Let $\xi\in\Xi$ and set $\xi'=\xi\vee\xih_n\in\Xi$. Monotonicity gives $F(\xi')\ge F(\xi)$, while
$|\xi'_k-\xih_{nk}|\le|\xi_k-\xih_{nk}|$ for each $k$, so the distance from $\xih_n$ does not increase. It is therefore enough to maximize over the upper box $B_n=[\xih_n,\xib]$. On $B_n$, the penalty equals
$\lambda\sum_kw_k(\xi_k-\xih_{nk})$, which is linear, so the objective is convex and continuous. A convex function
attains its maximum over a polytope at a vertex, and the vertices of $B_n$ are exactly the upper stress corners
$\xih_n^S$.
\end{proof}

\begin{remark}[Why supermodularity matters]\label[remark]{rem:hardness}
Convexity and monotonicity reduce the search to at most $2^m$ corners, but this can still be very expensive. In
fact, the inner maximization problem can be NP-hard even when $F$ is convex and nondecreasing. To see this, take a
graph with weighted Laplacian $L$ and maximum weighted degree $\Delta$, and let $\Xi=[0,1]^m$, $\xih_n=0$,
$F(x)=x^\top Lx+2\Delta\,\mathbf 1^\top x$, and $\lambda w_k=2\Delta$. Then $F$ is convex and nondecreasing on $\Xi$,
since $(Lx)_i\ge-\Delta$ there, and $F(x)-\lambda\|x-\xih_n\|_w=x^\top Lx$. The maximum of this function over
$\{0,1\}^m$ is the maximum cut \cite{Karp1972}. This $F$ is submodular rather than supermodular. Supermodularity
allows a further reduction through submodular minimization.

The weighted $\ell_1$ distance is important because it gives a linear penalty on the upper box
$[\xih_n,\xib]$. With a nonseparable distance such as $\ell_2$, the reduction to $\xi\ge\xih_n$ still holds, but the
objective on $[\xih_n,\xib]$ is a difference of convex functions and need not attain its maximum at a corner. The
weights $w_k$ allow different costs for changing different source intensities.
\end{remark}

\begin{proposition}[Supermodular stress function]\label[proposition]{prop:supermodular-set}
For every $n$ and $\lambda\ge0$, the function $S\mapsto\Psi_n(S;\lambda)$ is supermodular on $2^{[m]}$. Its
maximizers form a lattice under $\cup$ and $\cap$; in particular, there is a largest maximizer $\Smax_n(\lambda)$.
Both $\varphi_n(\lambda)$ and $\Smax_n(\lambda)$ can be computed by submodular minimization with a number of evaluations of $F$ (that is, PDE solves) polynomial in $m$
\cite{Schrijver2000,IwataFleischerFujishige2001,Orlin2009,Fujishige2005}.
\end{proposition}
\begin{proof}
Because $\xib\ge\xih_n$, the map from a set $S$ to its stress corner preserves unions and intersections:
$\xih_n^{S\cup T}=\xih_n^S\vee\xih_n^T$ and $\xih_n^{S\cap T}=\xih_n^S\wedge\xih_n^T$. Thus $S\mapsto F(\xih_n^S)$ is
supermodular, while $c_n$ is modular, so $\Psi_n$ is supermodular. If $A$ and $B$ are maximizers, then
\[
  \Psi_n(A\cup B)+\Psi_n(A\cap B)\ge\Psi_n(A)+\Psi_n(B)=2\varphi_n(\lambda),
\]
so $A\cup B$ and $A\cap B$ are also
maximizers. Therefore, the union of all maximizers is the largest maximizer. Finally, maximizing $\Psi_n$ is
equivalent to minimizing the submodular function $-\Psi_n(\cdot;\lambda)$, and the smallest and largest minimizers of
a submodular function can be computed in strongly polynomial time \cite{Fujishige2005}.
\end{proof}

Here \emph{an oracle call} is a submodular-minimization call that returns the largest minimizer of
$-\Psi_n(\cdot;\lambda)$, that is, the largest maximizer $\Smax_n(\lambda)$. All complexity statements below count
such calls.

\begin{theorem}[Nested stress chains]\label[theorem]{thm:chain}
For every $n$, if $0\le\lambda_1<\lambda_2$, then $\Smax_n(\lambda_2)\subseteq\Smax_n(\lambda_1)$. Thus, as $\lambda$
increases, the largest maximizers form a nested sequence of sets, the stress chain
$\mathcal C_n=\{\Smax_n(\lambda):\lambda\ge0\}$, which starts at $\Smax_n(0)=[m]$, ends at $K_n$, and has at most $m+1$ elements. Every $\Smax_n(\lambda)$ contains $K_n$. Moreover,
\[
  \varphi_n(\lambda)=\max_{S\in\mathcal C_n}\Psi_n(S;\lambda)\qquad(\lambda\ge0)
\]
is the maximum of at most $m+1$ affine functions, and $\varphi_n$ has at most $|\mathcal C_n|-1\le m$
breakpoints in $(0,\infty)$.
\end{theorem}
\begin{proof}
Let $A=\Smax_n(\lambda_1)$ and $B=\Smax_n(\lambda_2)$. Supermodularity at $\lambda_2$ and the maximality of $B$ give
$\Psi_n(A\cup B;\lambda_2)-\Psi_n(A;\lambda_2)\ge\Psi_n(B;\lambda_2)-\Psi_n(A\cap B;\lambda_2)\ge0$. Since $c_n$ is
nondecreasing,
\[
  \Psi_n(A\cup B;\lambda_1)-\Psi_n(A;\lambda_1)
  =\Psi_n(A\cup B;\lambda_2)-\Psi_n(A;\lambda_2)+(\lambda_2-\lambda_1)\bigl(c_n(A\cup B)-c_n(A)\bigr)\ge0,
\]
so $A\cup B$ is a maximizer at $\lambda_1$. Since $A$ is the largest maximizer, $A\cup B\subseteq A$; that is,
$B\subseteq A$. At $\lambda=0$, $[m]$ is a maximizer because $F$ is nondecreasing. For $\lambda$ larger than $\max\{(F(\xih_n^S)-F(\xih_n))/c_n(S):c_n(S)>0\}$, the maximizers are exactly the sets with $c_n(S)=0$, namely the
subsets of $K_n$, all of which give the same point $\xih_n$; the largest of them is $K_n$. A strictly decreasing
chain in $2^{[m]}$ has at most $m+1$ elements. For every $\lambda$, a set in the chain, namely $\Smax_n(\lambda)$,
attains the maximum $\varphi_n(\lambda)$, so it is enough to take the maximum over the chain. If $S$ is a maximizer,
then so is $S\cup K_n$, because adding indices of $K_n$ changes neither the stress corner nor its cost; hence every
$\Smax_n(\lambda)$ contains $K_n$. Distinct elements of $\mathcal C_n$ therefore differ only outside $K_n$ and have
distinct costs. Hence, for each pair of adjacent elements, there is exactly one value $\lambda\ge0$ at which the
largest maximizer changes from one to the other. If this transition value is positive, it is a breakpoint of
$\varphi_n$. If it is $0$, which happens when the two corners have equal damage, it need not give a breakpoint in
$(0,\infty)$: for example, if $m=1$, $\xih_{n1}<\xib_1$, and $F$ is constant, then $\mathcal C_n=\{[1],\emptyset\}$ but
$\varphi_n$ is constant. In either case, $\varphi_n$ has at most $|\mathcal C_n|-1$ breakpoints in $(0,\infty)$.
\end{proof}

\begin{remark}[Parametric submodular minimization and source ranking]\label[remark]{rem:ranking}
The nesting property in \cref{thm:chain} is a standard result of parametric submodular minimization
\cite{Topkis1978,Fujishige2005,GalloGrigoriadisTarjan1989}, and the proof above is Topkis's argument. Our PDE
results provide the supermodular function needed to use it. In practical terms, the chain $\mathcal C_n$ ranks the
sources for observation $n$. When $\lambda$ is small, changing source intensities carries a low penalty, so the
worst case raises every source to its upper bound. As $\lambda$ increases, sources leave the stressed set in a nested
order. Sources that remain in the set for larger $\lambda$ contribute more strongly to the worst-case objective relative to their transport costs. In \cref{sec:numerics}, the two sources aligned with the protected zone never leave before the other two, and in seven of the eight observations they are the last to leave.
\end{remark}

\begin{theorem}[Exact robust value and worst-case distribution]\label[theorem]{thm:exact}
Let $\rho\ge0$ and $\Lambda(\lambda)=\lambda\rho+\frac1N\sum_{n=1}^N\varphi_n(\lambda)$. Then
\begin{equation}\label{eq:exact}
  J(\rho)=\min_{\lambda\ge0}\Lambda(\lambda)=\min_{\lambda\in\mathcal B}\Lambda(\lambda),\qquad
  \mathcal B=\{0\}\cup\{\text{breakpoints of }\varphi_1,\dots,\varphi_N\},\quad|\mathcal B|\le Nm+1 .
\end{equation}
Let $\lambda^\star$ be a minimizer, and let $H_n=\Smax_n(\lambda^\star)$ and
$L_n=\lim_{\lambda\downarrow\lambda^\star}\Smax_n(\lambda)$ (the chain is eventually constant), so that $L_n\subseteq H_n$.
There are weights $p\in[0,1]^N$, at most one of them fractional, such that
\begin{equation}\label{eq:worstcase}
  P^\star=\frac1N\sum_{n=1}^N\Bigl(p_n\,\delta_{\xih_n^{H_n}}+(1-p_n)\,\delta_{\xih_n^{L_n}}\Bigr)
\end{equation}
is a worst-case distribution: $W(P^\star,\widehat P_N)\le\rho$ and $\mathbb E_{P^\star}F=J(\rho)$. In particular, for every
radius $\rho\ge0$, there is a worst-case distribution supported on the fixed set
$\mathcal C=\bigcup_n\{\xih_n^S:S\in\mathcal C_n\}$ of chain corners, which contains at most $N(m+1)$ points and does
not depend on $\rho$, and $P^\star$ uses at most $N+1$ of these points. Moreover, $J$ is the minimum of at most $Nm+1$
affine functions of $\rho$, with slopes in $\mathcal B$; hence it is concave and piecewise linear on $[0,\infty)$.
\end{theorem}

\begin{proof}
\emph{Upper bound.} Let $W(P,\widehat P_N)\le\rho$, $\varepsilon>0$, and $\gamma\in\Gamma(P,\widehat P_N)$ with
$\int\|\xi-\zeta\|_w\,\mathrm d\gamma\le\rho+\varepsilon$. For $\lambda\ge0$, by \cref{thm:corners},
\[
  \mathbb E_PF=\int\bigl[F(\xi)-\lambda\|\xi-\zeta\|_w\bigr]\mathrm d\gamma+\lambda\int\|\xi-\zeta\|_w\,\mathrm d\gamma
  \le\frac1N\sum_n\varphi_n(\lambda)+\lambda(\rho+\varepsilon).
\]
Letting $\varepsilon\to0$ gives $J(\rho)\le\Lambda(\lambda)$ for every $\lambda\ge0$.

\emph{Location of a minimizer.} By \cref{thm:chain}, $\Lambda$ is convex and piecewise linear with kinks only in
$\mathcal B\setminus\{0\}$. Beyond the last breakpoint its slope is $\rho-\frac1N\sum_nc_n(K_n)=\rho\ge0$, so
$\Lambda$ attains its minimum over $[0,\infty)$ at a point of $\mathcal B$.

\emph{One-sided slopes.} On a small interval $(\lambda^\star,\lambda^\star+\varepsilon)$ the largest maximizer is
constant, equal to $L_n$, so $\varphi_n=\Psi_n(L_n;\cdot)$ there; by continuity, $L_n$ is a maximizer at
$\lambda^\star$ and $\varphi_n'(\lambda^\star{+})=-c_n(L_n)$. If $\lambda^\star>0$, the same argument on the left
yields a set $H'$ that is a maximizer at $\lambda^\star$, hence $H'\subseteq H_n$, while
$H_n\subseteq H'$ by \cref{thm:chain}; thus $\varphi_n'(\lambda^\star{-})=-c_n(H_n)$, and $H_n$ is a maximizer.

\emph{Construction of a worst-case distribution.} If $\lambda^\star>0$, optimality gives
$\Lambda'(\lambda^\star{-})\le0\le\Lambda'(\lambda^\star{+})$, that is,
\[
  \frac1N\sum_nc_n(L_n)\le\rho\le\frac1N\sum_nc_n(H_n).
\]
We choose $p_n$ successively to use the available transportation budget: with $\Delta_n=c_n(H_n)-c_n(L_n)\ge0$, choose $p_1,p_2,\dots$ in turn as large as
possible in $[0,1]$ subject to $\frac1N\sum_n[c_n(L_n)+p_n\Delta_n]\le\rho$. Equality is reached, and at most one $p_n$
is fractional. Coupling each $\xih_n$ with the two corners in~\eqref{eq:worstcase} costs exactly $\rho$, so
$W(P^\star,\widehat P_N)\le\rho$. Since $H_n$ and $L_n$ both maximize $\Psi_n(\cdot;\lambda^\star)$, we have
$F(\xih_n^{H_n})=\varphi_n(\lambda^\star)+\lambda^\star c_n(H_n)$, and similarly for $L_n$. It follows that
\[
  \mathbb E_{P^\star}F=\frac1N\sum_n\varphi_n(\lambda^\star)+\lambda^\star\rho=\Lambda(\lambda^\star)\ge J(\rho).
\]
If $\lambda^\star=0$, we take $p=0$: the cost is $\frac1N\sum_nc_n(L_n)\le\rho$ (from $\Lambda'(0{+})\ge0$) and the
value is $\frac1N\sum_n\varphi_n(0)=\Lambda(0)$. In both cases $J(\rho)=\Lambda(\lambda^\star)$. Finally, every observation with $p_n\in\{0,1\}$ contributes one support
point and the split observation, if any, contributes two, so $P^\star$ has at most $N+1$ support points. Since
$H_n=\Smax_n(\lambda^\star)$ and $L_n$ belong to $\mathcal C_n$, all support points lie in $\mathcal C$. Since $\mathcal B$ does not depend on $\rho$, \eqref{eq:exact} expresses $J$ as the minimum of the affine functions
$\rho\mapsto\lambda\rho+\frac1N\sum_n\varphi_n(\lambda)$, $\lambda\in\mathcal B$.
\end{proof}

\begin{remark}[What supermodularity adds to the support result]\label[remark]{rem:novelty}
For a Wasserstein ball around an empirical distribution, a worst-case distribution with at most $N+1$ support
points, in which at most one observation is split, exists under general conditions \cite{GaoKleywegt2023}; we do
not claim this cardinality bound as new. What \cref{thm:exact} adds under monotonicity, convexity, and
supermodularity is information about the location of the support points and the form of the dual function.
(i) The candidate support set $\mathcal C$ has at most $N(m+1)$ points, is computed once by \cref{alg:chain}, and
serves every radius. Under monotonicity and convexity alone, the corresponding radius-independent set consists of
all $N2^m$ upper stress corners (\cref{thm:corners}); for a general loss, the support points are maximizers of
$F(\xi)-\lambda^\star\|\xi-\xih_n\|_w$ over $\Xi$ and may change with $\rho$ through $\lambda^\star$. (ii) The two points of
a split observation are adjacent elements $H_n\supseteq L_n$ of the same chain, so the worst case raises the sources
in $H_n\setminus L_n$ to their upper bounds with probability $p_n$. (iii) $J$ is available for all radii at once as
the minimum of at most $Nm+1$ explicitly known affine functions, as in \cref{fig:dro} (right). (iv) Without
supermodularity, even evaluating $\varphi_n(\lambda)$ for a convex nondecreasing $F$ can be NP-hard
(\cref{rem:hardness}).
\end{remark}

\begin{remark}[Weighted observations and observation-dependent costs]\label[remark]{rem:weighted}
The results of this section hold, with the same proofs, in two more general settings. First, $\widehat P_N$ may be
replaced by $\sum_n\alpha_n\delta_{\xih_n}$ with $\alpha_n>0$ and $\sum_n\alpha_n=1$, for example to merge repeated
observations or to weight observations by reliability. The factor $\frac1N$ is then replaced by $\alpha_n$ in
$\Lambda$, in~\eqref{eq:worstcase}, and in the transportation budget used to choose the $p_n$, and the greedy choice
still leaves at most one $p_n$ fractional. Second, for distinct observations, the cost of moving mass away from
$\xih_n$ may use observation-dependent weights $w^{(n)}>0$, that is, $\|\xi-\xih_n\|_{w^{(n)}}$. Then
$c_n(S)=\sum_{k\in S}w^{(n)}_k(\xib_k-\xih_{nk})$ is still modular and nondecreasing, which is all that
\cref{thm:corners,prop:supermodular-set,thm:chain} use. In this case $W$ is an optimal transport cost rather than a
metric, and \cref{thm:exact} holds with $W$ replaced by this cost, because its proof uses only the coupling bound
and the explicit construction.
\end{remark}

\begin{remark}[Degenerate cases]\label[remark]{rem:degenerate}
The proof of \cref{thm:exact} covers the following cases without change. (i) If $\rho=0$, then $\Lambda$ is
nonincreasing and is constant beyond the largest point of $\mathcal B$, where it attains its minimum; there
$L_n=K_n$, the budget forces $p_n=0$ whenever $c_n(H_n)>c_n(L_n)$, and $P^\star=\widehat P_N$, so
$J(0)=\frac1N\sum_nF(\xih_n)$. (ii) If $\Lambda$ attains its minimum at several points of $\mathcal B$, the
construction may use any of them. (iii) Breakpoints of different $\varphi_n$ may coincide; then several observations
have $c_n(H_n)>c_n(L_n)$, and the greedy choice still leaves at most one $p_n$ fractional. (iv) When the maximizer of $\Psi_n(\cdot;\lambda)$ is not unique, the construction uses the largest maximizer and its right limit, both of which
belong to $\mathcal C_n$; sets that differ only in $K_n$ give the same corner. \Cref{sec:numerics} tests all four cases
against a linear program.
\end{remark}

\begin{algorithm}[t]
\caption{Stress chain and breakpoints of $\varphi_n$ by intersection search}\label{alg:chain}
\begin{algorithmic}[1]
\Require oracle $\mathcal O_n(\lambda)=\Smax_n(\lambda)$ (one oracle call); evaluations of $F$
\State $\mathcal C_n\gets\{[m],K_n\}$, \ $\mathcal B_n\gets\emptyset$; \textbf{if} $[m]=K_n$ \textbf{return}
\State \Call{Refine}{$[m],K_n$}; \Return $(\mathcal C_n,\mathcal B_n)$
\Procedure{Refine}{$S_a,S_b$}\Comment{$S_a\supsetneq S_b$: adjacent known chain sets}
  \State $\lambda_c\gets\dfrac{F(\xih_n^{S_a})-F(\xih_n^{S_b})}{c_n(S_a)-c_n(S_b)}$\Comment{intersection of the two lines}
  \State $S_c\gets\mathcal O_n(\lambda_c)$
  \If{$\Psi_n(S_c;\lambda_c)>\Psi_n(S_a;\lambda_c)$}
     \State $\mathcal C_n\gets\mathcal C_n\cup\{S_c\}$; \Call{Refine}{$S_a,S_c$}; \Call{Refine}{$S_c,S_b$}
  \Else
     \State $\mathcal B_n\gets\mathcal B_n\cup\{\lambda_c\}$\Comment{transition value; a breakpoint if $\lambda_c>0$}
  \EndIf
\EndProcedure
\end{algorithmic}
\end{algorithm}

\begin{proposition}[Oracle complexity]\label[proposition]{prop:algorithm}
\Cref{alg:chain} returns $\mathcal C_n$ and the set $\mathcal B_n$ of its transition values, which consists of the breakpoints of $\varphi_n$ and possibly $0$, using at most
$2|\mathcal C_n|-3\le2m-1$ oracle calls, and none if $|\mathcal C_n|=1$. Hence $J(\rho)$ and a worst-case
distribution are obtained from at most $N(2m-1)$ submodular minimizations, followed by the evaluation
of~\eqref{eq:exact} at the $|\mathcal B|\le Nm+1$ points of $\mathcal B=\{0\}\cup\bigcup_n\mathcal B_n$.
\end{proposition}
\begin{proof}
We use the following property at each recursive step: $S_a=\Smax_n(\lambda_a)$ and $S_b=\Smax_n(\lambda_b)$ for some
$\lambda_a<\lambda_b$, and $c_n(S_a)>c_n(S_b)$, since $S_a\setminus S_b$ is nonempty and disjoint from $K_n$. Let $\Psi_a=\Psi_n(S_a;\cdot)$ and $\Psi_b=\Psi_n(S_b;\cdot)$; these functions are affine. The difference $\Psi_a-\Psi_b$ is
strictly decreasing, nonnegative at $\lambda_a$, and nonpositive at $\lambda_b$, so the intersection point $\lambda_c$
lies in $[\lambda_a,\lambda_b]$. If $\varphi_n(\lambda_c)=\Psi_a(\lambda_c)$, then the convex function $\varphi_n\ge\Psi_a$ coincides with $\Psi_a$ at $\lambda_a$ and $\lambda_c$, hence on $[\lambda_a,\lambda_c]$; likewise, it coincides with $\Psi_b$
on $[\lambda_c,\lambda_b]$. In this case, $\lambda_c$ is the only transition value in $[\lambda_a,\lambda_b)$, and no chain set lies
strictly between $S_a$ and $S_b$. Otherwise, the oracle returns a new chain set $S_c$ with
$S_b\subsetneq S_c\subsetneq S_a$ (\cref{thm:chain}) and $\lambda_a<\lambda_c<\lambda_b$, and the property holds for both
recursive calls. Each oracle call therefore either finds one of the $|\mathcal C_n|-2$ interior chain sets or
confirms one of the $|\mathcal C_n|-1$ transition values.
\end{proof}

This bound counts submodular-minimization calls, not PDE solves. Each evaluation of $F$ requires a PDE solve,
and the number of evaluations per call depends on the submodular-minimization algorithm. \Cref{tab:scaling}
reports both counts.

\begin{remark}[The role of concavity of $g$]
If $g$ is convex and $\Phi$ is linear, then $F$ is concave (\cref{cor:param}), and the inner maximization problems of \cref{thm:corners} need not attain their maximum at a corner. Together with \cref{thm:convex-reaction}, this shows that the sign of
$g''$ determines whether the structure used in this section is available.
\end{remark}

\section{Application: robust remediation capacity with Monod biodegradation}\label{sec:application}

Let $\Omega\subset\R^2$ be a bounded Lipschitz domain. We model the contaminant concentration $u$ by
\begin{equation}\label{eq:remediation}
  -\nabla\!\cdot(D\nabla u)+\beta\cdot\nabla u+\gamma u+g_M(u)=\sum_{k=1}^m\xi_kf_k-\sum_{j=1}^Jy_jr_j\ \text{ in }\Omega,
  \qquad u=0\ \text{ on }\partial\Omega.
\end{equation}
Here $D$ is uniformly elliptic, $\beta\in W^{1,\infty}$ is an incompressible velocity field ($\operatorname{div}\beta=0$),
and $\gamma\ge0$. The nonnegative functions $f_k$ represent uncertain contaminant sources. Each nonnegative function
$r_j$, with $\int r_j=1$, describes the spatial footprint of treatment station $j$, and $y_j$ is its removal
capacity. The capacities $y$ belong to a convex compact set $Y\subset\R^J_+$. We use the Monod biodegradation
reaction \cite{Monod1949}
\[
  g_M(r)=\frac{\kappa r}{K_m+r}\ (r\ge0),\qquad g_M(r)=\frac{\kappa}{K_m}\,r\ (r<0).
\]
This function is nondecreasing, concave, and $\kappa/K_m$-Lipschitz on $\R$, as required in \cref{sec:order}.
Thus the Monod reaction itself satisfies the concavity condition; it is not an additional assumption made for the proofs. For negative arguments,
we use a linear extension only to meet the assumptions of the general theory in \cref{sec:framework}; it is not intended to model biodegradation at negative concentrations. \Cref{prop:positivity} and \cref{sec:numerics} address nonnegativity of the states that determine the robust risk. We measure the damage in a protected zone $Z$ by
\[
  \Phi(x,r)=\mathbf 1_Z(x)(r-u_{\rm safe})_+^2
\]
and write
the resulting damage as $F(\xi,y)$.

\begin{proposition}[Convex robust design]\label[proposition]{prop:design}
\begin{enumerate}[label=(\alph*)]
\item $F$ is jointly convex in $(\xi,y)$, nondecreasing in $\xi$, nonincreasing in $y$, and supermodular in
  $(\xi,-y)$.
\item For each $y$, $F(\cdot,y)$ satisfies the assumptions of \cref{sec:dro}, so the robust risk
  $R_\rho(y)=\sup\{\mathbb E_PF(\xi,y):W(P,\widehat P_N)\le\rho\}$ is computed exactly by \cref{thm:exact}.
\item $R_\rho$ is convex on $Y$, so $\min_{y\in Y}R_\rho(y)$ is a convex program. If $P^\star_y$ is a worst-case distribution at $y$, for example the one from \cref{thm:exact}, then any subgradient of $y'\mapsto\mathbb E_{P^\star_y}F(\xi,y')$ at $y$
  is a subgradient of $R_\rho$ at $y$.
\end{enumerate}
\end{proposition}
\begin{proof}
Equation~\eqref{eq:remediation} is covered by \cref{ex:local}, and its load $\ell(\xi,y)$ depends affinely on $\xi$
and $y$. (a) Joint convexity follows from \cref{thm:order}\ref{it:convex} and the argument of \cref{thm:damage}. To
see monotonicity and supermodularity, write $z=-y$. In the variables $(\xi,z)$, all load coefficients $f_k$ and
$r_j$ are nonnegative, so \cref{cor:param,thm:damage} apply. (b) This is part~(a) with $y$ fixed; the term
$-\sum_jy_jr_j$ is then part of $\ell_0$. (c) For every fixed $P$, the expected damage $y\mapsto\mathbb E_PF(\xi,y)$ is
convex, and the uncertainty set does not depend on $y$. A supremum of convex functions is convex. If $\sigma$ is a subgradient of $y'\mapsto\mathbb E_{P^\star_y}F(\xi,y')$ at $y$, then $R_\rho(y')\ge\mathbb E_{P^\star_y}F(\xi,y')\ge R_\rho(y)+\ip{\sigma}{y'-y}$.
\end{proof}

\paragraph{Differentiability and adjoint gradients}
The cutting-plane method needs subgradients of $R_\rho$, and by \cref{prop:design}(c) these come from derivatives of
$F(\xi,\cdot)$. We prove that $F$ is continuously differentiable and derive its gradient using an adjoint equation. Let $V=H^1_0(\Omega)$, let $\E$ be the form of \cref{ex:local} for~\eqref{eq:remediation}, and write
$\ell(\zeta,\eta)=\sum_k\zeta_kf_k-\sum_j\eta_jr_j\in V^*$ for $(\zeta,\eta)\in\R^m\times\R^J$.

\begin{lemma}[Differentiability of the state and the damage]\label[lemma]{lem:diff}
Let $f_k,r_j\in V^*$. The solution map $\mathcal T:(\xi,y)\mapsto u_{\xi,y}$ of~\eqref{eq:remediation} is continuously Fr\'echet
differentiable from $\R^m\times\R^J$ to $V$. Its derivative at $(\xi,y)$ in the direction $(\zeta,\eta)$ is the
unique $v\in V$ with
\begin{equation}\label{eq:linearized}
  \E(v,\varphi)+\int_\Omega g_M'(u)\,v\,\varphi\dx=\ip{\ell(\zeta,\eta)}{\varphi}\qquad\forall\varphi\in V,\qquad u=u_{\xi,y}.
\end{equation}
The damage $F$ is also continuously differentiable on $\R^m\times\R^J$, and $\partial F/\partial y_j=-\ip{r_j}{\psi}$, where $\psi\in V$ uniquely solves the adjoint problem
\begin{equation}\label{eq:adjoint}
  \E(v,\psi)+\int_\Omega g_M'(u)\,v\,\psi\dx=\int_Z2\bigl(u-u_{\rm safe}\bigr)_+v\dx\qquad\forall v\in V .
\end{equation}
\end{lemma}
\begin{proof}
\emph{Properties of $g_M$.} The two branches of $g_M$ have the same value and the same slope $\kappa/K_m$ at $r=0$, so
$g_M\in C^1(\R)$ with $0\le g_M'\le\kappa/K_m$. For $r>0$, $|g_M''(r)|=2\kappa K_m/(K_m+r)^3\le2\kappa/K_m^2=:L'$, and
$g_M''=0$ for $r<0$. Hence $g_M'$ is Lipschitz with constant $L'$, and
\begin{equation}\label{eq:taylor}
  \bigl|g_M(r+s)-g_M(r)-g_M'(r)s\bigr|\le\tfrac{L'}2\,s^2\qquad(r,s\in\R).
\end{equation}
Since $\Omega\subset\R^2$ is bounded, $V$ embeds continuously in $L^4(\Omega)$; let $\|v\|_{L^4}\le C_4\|v\|_V$.

\emph{Linearized operators.} For $w\in V$, define $B_w:V\to V^*$ by
$\ip{B_wv}{\varphi}=\E(v,\varphi)+\int_\Omega g_M'(w)v\varphi\dx$. Because $0\le g_M'(w)\le\kappa/K_m$, the operator $B_w$
is bounded and $\ip{B_wv}{v}\ge\alpha\|v\|_V^2$. By the Lax--Milgram lemma, $B_w$ is invertible with
$\|B_w^{-1}\|\le\alpha^{-1}$ for every $w\in V$, and the same holds for its adjoint. In particular,
\eqref{eq:linearized} and~\eqref{eq:adjoint} have unique solutions.

\emph{Differentiability of $\mathcal T$.} Fix $(\xi,y)$ and an increment $(\zeta,\eta)$, and let $u=\mathcal T(\xi,y)$,
$\tilde u=\mathcal T(\xi+\zeta,y+\eta)$, $\delta=\tilde u-u$, and $v=B_u^{-1}\ell(\zeta,\eta)$. By \cref{lem:wellposed},
$\|\delta\|_V\le\alpha^{-1}\|\ell(\zeta,\eta)\|_{V^*}\le C|(\zeta,\eta)|$. Subtracting the two state equations gives
\[
  \ip{B_u(\delta-v)}{\varphi}=-\int_\Omega q\,\varphi\dx\quad\forall\varphi\in V,\qquad
  q=g_M(\tilde u)-g_M(u)-g_M'(u)\,\delta .
\]
By~\eqref{eq:taylor} and H\"older's inequality,
$|\int q\varphi|\le\frac{L'}2\|\delta\|_{L^4}^2\|\varphi\|_{L^2}\le\frac{L'}2C_4^2C_V\|\delta\|_V^2\|\varphi\|_V$. Therefore
$\|\delta-v\|_V\le\frac{L'}{2\alpha}C_4^2C_V\|\delta\|_V^2=O(|(\zeta,\eta)|^2)$, so $\mathcal T$ is Fr\'echet differentiable with
$D\mathcal T(\xi,y)(\zeta,\eta)=B_u^{-1}\ell(\zeta,\eta)$. For $u,w\in V$, H\"older's inequality also gives
$\|B_u-B_w\|_{\mathcal L(V,V^*)}\le L'C_VC_4^2\|u-w\|_V$, and
$\|B_u^{-1}-B_w^{-1}\|=\|B_u^{-1}(B_w-B_u)B_w^{-1}\|\le\alpha^{-2}\|B_u-B_w\|$. Since $\mathcal T$ is Lipschitz, $D\mathcal T$ is
continuous.

\emph{The damage and the adjoint.} Let $\Phi_0(r)=(r-u_{\rm safe})_+^2$. Its derivative $\Phi_0'(r)=2(r-u_{\rm safe})_+$ is
Lipschitz with constant $2$, so $|\Phi_0(r+s)-\Phi_0(r)-\Phi_0'(r)s|\le s^2$. Hence $u\mapsto\int_Z\Phi_0(u)\dx$ is
continuously Fr\'echet differentiable on $L^2(\Omega)$ with derivative $v\mapsto\int_Z\Phi_0'(u)v\dx$, and by the chain
rule and $V\hookrightarrow L^2$, $F$ is continuously differentiable. For the direction $(0,e_j)$, let
$v_j=B_u^{-1}\ell(0,e_j)$. Testing~\eqref{eq:adjoint} with $v=v_j$ and~\eqref{eq:linearized} with $\varphi=\psi$ gives
\[
  \frac{\partial F}{\partial y_j}=\int_Z\Phi_0'(u)\,v_j\dx=\E(v_j,\psi)+\int_\Omega g_M'(u)\,v_j\,\psi\dx
  =\ip{\ell(0,e_j)}{\psi}=-\ip{r_j}{\psi}.
\]
\end{proof}

The argument only requires the embedding $H^1_0(\Omega)\hookrightarrow L^4(\Omega)$ and the Lipschitz continuity of $g_M'$ and
$\Phi_0'$, so \cref{lem:diff} holds for bounded domains in $\R^d$ with $d\le4$. For the finite-difference model of
\cref{sec:numerics}, the same conclusions follow from the implicit function theorem in $\R^n$, and $\psi$ is computed
with the transposed Jacobian. When $r_j\in L^2(\Omega)$, as in our experiments, $\ip{r_j}{\psi}=\int_\Omega r_j\psi\dx$.
For a worst-case distribution $P^\star_y=\sum_l\omega_l\delta_{\xi^l}$ from
\cref{thm:exact}, with adjoint states $\psi^l$, we combine the adjoint gradients of its support points.
\Cref{prop:design}(c) then gives the robust subgradient
\begin{equation}\label{eq:robustgrad}
  \sigma_j(y)=-\sum_l\omega_l\int_\Omega r_j\,\psi^l\dx ,\qquad j=1,\dots,J.
\end{equation}
Thus at most $N+1$ adjoint solves are needed, one per support point.

\begin{proposition}[Complementarity and diminishing returns]\label[proposition]{prop:complements}
\begin{enumerate}[label=(\alph*)]
\item For all $j,k$ and $t>0$, the damage reduction $F(\xi,y)-F(\xi,y+te_j)$ is nondecreasing in $\xi_k$. Thus, a
  given increase in treatment capacity reduces more damage when the contaminant load is higher.
\item $F(\xi,\cdot)$ is supermodular. For capacities $\bar y_j>0$ and $y^T=\sum_{j\in T}\bar y_je_j$, the siting benefit
  $B_\xi(T)=F(\xi,0)-F(\xi,y^T)$ is nondecreasing and submodular on $2^{[J]}$, and so is $\mathbb E_PB_\xi(T)$ for any
  fixed distribution $P$. Therefore, under a constraint on the number of selected stations, the greedy method
  achieves at least a fraction $1-1/e$ of the optimal nominal or expected benefit \cite{NemhauserWolseyFisher1978}.
\end{enumerate}
\end{proposition}
\begin{proof}
(a) This is the increasing-differences property of $(\xi,z)\mapsto F(\xi,-z)$. (b) Supermodularity is preserved
under $y\mapsto-y$, which exchanges $\vee$ and $\wedge$. Since $T\mapsto y^T$ preserves unions and intersections,
$B_\xi$ is submodular; it is nondecreasing because $F$ is nonincreasing in $y$. Expectations preserve both
properties.
\end{proof}

\begin{remark}
The worst-case siting benefit $\min_P\mathbb E_PB_\xi(T)$ is a minimum of submodular functions and need not
itself be submodular \cite{KrauseEtAl2008}. We therefore optimize continuous treatment capacities rather than
claim a greedy guarantee for the worst-case siting benefit.
\end{remark}

\paragraph{Nonnegative concentrations}
The concentration model~\eqref{eq:remediation} is physically meaningful only for nonnegative states, and nonnegativity cannot hold throughout $\Xi$:
if $\xil=0$ and $y\in Y\setminus\{0\}$, the load at $\xi=0$ is $-\sum_jy_jr_j\le0$ and nonzero, so
\cref{thm:order}\ref{it:mono} gives $u_{0,y}\le0$ with $u_{0,y}\ne0$. The robust risk, however, depends only on states at the upper stress corners (\cref{thm:corners,thm:exact}). We therefore use~\eqref{eq:remediation} as a surrogate model for risk optimization: the nonnegativity
guarantee below covers all scenarios needed to compute the worst-case objective, including the support of the
worst-case distribution in \cref{thm:exact}, but not every realization in $\Xi$. The following proposition gives a computable sufficient condition for nonnegativity throughout the upper boxes $[\xih_n,\xib]$. Let $u^n=u_{\xih_n,0}$, and let $v^n_j\in V$ be the
derivative of $u_{\xi,y}$ with respect to $y_j$ at $(\xih_n,0)$, that is, the solution of~\eqref{eq:linearized} with
$u=u^n$ and load $-r_j$.

\begin{proposition}[Nonnegative states in the robust problem]\label[proposition]{prop:positivity}
Let $Y_+=\{y\in Y:\ u^n+\sum_jy_jv^n_j\ge0\ \text{a.e.\ in }\Omega,\ n=1,\dots,N\}$. Then $Y_+$ is convex and closed,
and $u_{\xi,y}\ge0$ a.e.\ for every $y\in Y_+$ and every $\xi\in\bigcup_n[\xih_n,\xib]$. Consequently,
$\min_{y\in Y_+}R_\rho(y)$ is a convex program in which every state that determines $R_\rho$ is nonnegative.
\end{proposition}
\begin{proof}
The constraints are affine in $y$, so $Y_+$ is convex and closed. Let $y\in Y_+$ and $\xi\in[\xih_n,\xib]$. Since
$\ell(\xi,y)-\ell(\xih_n,y)=\sum_k(\xi_k-\xih_{nk})f_k\in\Vsp$, \cref{thm:order}\ref{it:mono} gives
$u_{\xi,y}\ge u_{\xih_n,y}$. The load is affine in $y$, so $t\mapsto u_{\xih_n,ty}$ is convex in the a.e.\ order by
\cref{thm:order}\ref{it:convex}, and hence $u_{\xih_n,y}-u^n\ge(u_{\xih_n,ty}-u^n)/t$ for $t\in(0,1]$. By
\cref{lem:diff}, the right-hand side converges in $V$, and hence in $L^2$, to $\sum_jy_jv^n_j$ as $t\downarrow0$, and
the cone of a.e.\ nonnegative functions is closed in $L^2$. Therefore $u_{\xi,y}\ge u^n+\sum_jy_jv^n_j\ge0$. Finally,
$R_\rho$ is convex by \cref{prop:design}(c), and by \cref{thm:corners,thm:exact} it depends on $F(\cdot,y)$ only
through its values at the upper stress corners, which lie in $\bigcup_n[\xih_n,\xib]$.
\end{proof}

To check whether $y\in Y_+$, we solve $N$ state equations and $NJ$ linearized equations. Since
$\sum_jy_jv^n_j\ge(\sum_jy_j)\min_jv^n_j$ for $y\ge0$, the single test $u^n+(\max_{y\in Y}\sum_jy_j)\min_jv^n_j\ge0$
for all $n$ shows that $Y_+=Y$. We use both tests in \cref{sec:numerics}.

\section{Numerical experiments}\label{sec:numerics}

These synthetic experiments test the theory and algorithms on a model problem; they are not a field study. The code
and data are available from the corresponding author upon reasonable request.

\paragraph{Setup} We discretize $\Omega=(0,1)^2$ on a uniform grid with $h=1/64$, using the five-point scheme for
diffusion, first-order upwind differences for advection, and the Monod reaction at the grid nodes, with $D=0.02$,
$\beta=(0.8,0)$, $\gamma=0.1$, $\kappa=2$, and $K_m=0.5$. At the nodes $x_i$ we solve $K_h^{\rm FD}u+g_M(u)=f$, where
$K_h^{\rm FD}$ is the finite-difference matrix and $f$ holds the nodal load values; this discretization
of~\eqref{eq:remediation} is consistent and first order in $h$. The discrete damage is $F_h=h^2\sum_i\Phi(x_i,u_i)$.
This is \cref{ex:discrete} with $\mu_i=h^2$ and $K=h^2K_h^{\rm FD}$, a Z-matrix with positive definite symmetric part
(\cref{app:examples}), so every theorem holds \emph{exactly} for $F_h$. We write $J_h(\rho)$ and $R_{\rho,h}$ for the
robust value and robust risk of the discrete model; we compute them exactly up to solver tolerances, and they differ
from $J(\rho)$ and $R_\rho$ only through the discretization, which we examine at the end of this section. We use $m=4$
Gaussian sources (amplitude $25$, width $0.05$) at $(0.15,0.20)$, $(0.22,0.42)$, $(0.15,0.62)$, and $(0.22,0.82)$, and
three candidate treatment stations at $(0.55,x_2)$, $x_2\in\{0.3,0.5,0.7\}$, with unit-mass footprints of width $0.1$.
The protected zone is $Z=[0.75,0.95]\times[0.35,0.65]$ with $u_{\rm safe}=0.35$, the uncertainty set is $\Xi=[0,2]^4$,
and $w=\mathbf 1$. The $N=8$ observations are $\xih_{nk}=\min\{0.9\,e^{0.35Z_{nk}},2\}$ with independent standard
normal $Z_{nk}$ (NumPy's default generator, seed $7$; see \cref{tab:obs}). Newton's method is run to a relative
residual of $10^{-12}$.

\paragraph{Testing supermodularity and convexity} For $150$ random pairs $\xi,\eta\in\Xi$, each with a random
$\theta\sim U[0,1]$, \cref{tab:structure} reports the smallest supermodularity gap
$F_h(\xi\vee\eta)+F_h(\xi\wedge\eta)-F_h(\xi)-F_h(\eta)$, the corresponding pointwise gap for $u$, and the smallest
convexity gap. For the Monod reaction, all gaps are nonnegative up to rounding error, as \cref{thm:order,thm:damage}
predict. For the strictly convex softplus reaction (\cref{ex:softplus}, $\kappa=2$) with linear damage, $91\%$ of the
pairs violate supermodularity. With the damage $(r-u_{\rm safe})_+^2$, no tested pair violates supermodularity of the
damage function, although the state map violates it; thus a failure for the state map need not cause a violation for
every damage function. \Cref{cor:tradeoff} gives a necessary curvature condition for this behavior in the
large-domain limit; whether such a condition is sufficient on general domains is open (\cref{sec:open}). For
identical sources $f_1=f_2$ and linear damage, \eqref{eq:strict} holds at all four tested base points, as
\cref{thm:convex-reaction} predicts: the second differences lie between $-9.2\times10^{-3}$ and $-3.0\times10^{-4}$,
compared with $+9.4\times10^{-4}$ to $+3.3\times10^{-2}$ for the Monod reaction.

\begin{table}[!t]
\centering
\caption{Minimum gaps over $150$ random pairs. Negative gaps larger than the numerical tolerance indicate a
violation of the stated property; values of order $10^{-16}$ are rounding errors. $\max F_h$, the largest damage
observed, indicates the scale.\label{tab:structure}}
\footnotesize\setlength{\tabcolsep}{4pt}
\begin{tabular}{llccccc}
\toprule
reaction & damage $\Phi$ & \shortstack{supermod.\\gap of $F_h$} & \shortstack{pointwise\\gap of $u$} & \shortstack{convexity\\gap of $F_h$} & $\max F_h$ & \shortstack{violating\\pairs}\\
\midrule
Monod (concave) & $(r-u_{\rm safe})_+^2$ & $-1.4\times10^{-17}$ & $-8.9\times10^{-16}$ & $6.5\times10^{-7}$ & $0.209$ & $0\%$\\
Monod (concave) & $r$ & $-6.9\times10^{-18}$ & $-8.9\times10^{-16}$ & $2.9\times10^{-6}$ & $0.130$ & $0\%$\\
softplus (convex) & $r$ & $-3.0\times10^{-3}$ & $-7.0\times10^{-2}$ & $-2.2\times10^{-3}$ & $0.070$ & $91\%$\\
softplus (convex) & $(r-u_{\rm safe})_+^2$ & $-3.5\times10^{-18}$ & $-8.8\times10^{-2}$ & $-1.2\times10^{-4}$ & $0.041$ & $0\%$\\
\bottomrule
\end{tabular}
\end{table}

\begin{table}[!t]
\centering
\caption{Stress chains and validation ($N=8$, $m=4$).\label{tab:dro}}
\small
\begin{tabular}{lc}
\toprule
chain lengths $|\mathcal C_n|$, $n=1,\dots,8$ & $4,3,3,4,2,3,4,4$\\
oracle calls of \cref{alg:chain} (bound $2|\mathcal C_n|-3$) & $5,3,3,5,1,3,5,5$\\
stress-chain corners $\sum_n|\mathcal C_n|$ vs.\ all corners $N2^m$ & $27$ vs.\ $128$\\
breakpoints of $\Lambda$ & $19$\\
$\max|J_h-J_h^{\rm LP}|$, corners only / corners $+$ $5^4$ grid & $5.6\times10^{-17}$ / $5.6\times10^{-17}$\\
$J_h(0.75)$, $\lambda^\star$, split observations & $0.0529$, $0.0666$, one ($p_4=0.548$)\\
$J_h(1.5)$, $\lambda^\star$, split observations & $0.0994$, $0.0542$, one ($p_6=0.065$)\\
distinct PDE solves: $121$ corners${}+625$ grid points${}-1$ shared & $745$\\
\bottomrule
\end{tabular}
\end{table}

\begin{figure}[!t]
\centering
\includegraphics[width=0.86\textwidth]{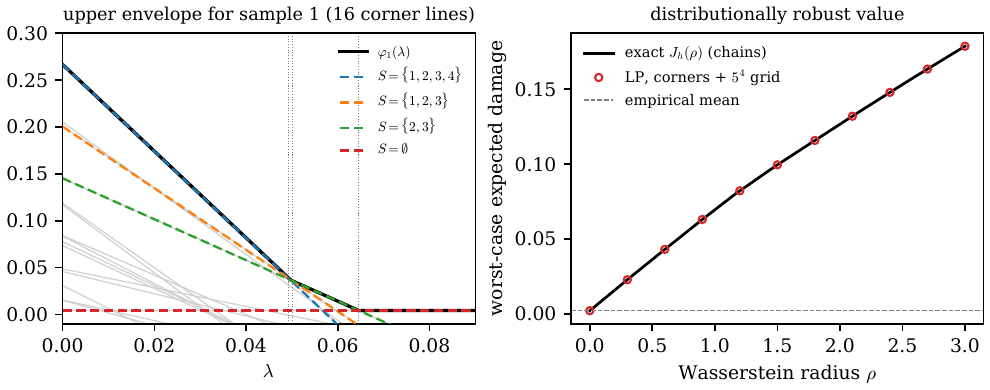}
\caption{Left: the $16$ corner lines $\lambda\mapsto\Psi_1(S;\lambda)$ (gray), the stress chain (dashed), and the upper
envelope $\varphi_1$ with its three breakpoints. Right: the robust value $J_h(\rho)$ of the discrete model from \cref{thm:exact}, compared with an independent linear-program calculation over $128+625$ candidate support points.\label{fig:dro}}
\end{figure}

\begin{table}[!t]
\centering
\caption{Robust remediation design by cutting planes (\cref{prop:design}(c)). The total capacity is $0.1$,
``cap'' is the capacity limit per site, and $y^\star$ lists the nonzero capacities by site location. In every run, the
master gap is at the level of rounding error ($|\mathrm{UB}-\mathrm{LB}|\le2\times10^{-16}\,\mathrm{UB}$); overall
accuracy is limited by the PDE and adjoint solves.\label{tab:solver}}
\small
\begin{tabular}{cccccp{2.35in}c}
\toprule
sites & $\rho$ & cap & $R_{\rho,h}$(uniform) & $R_{\rho,h}(y^\star)$ & $y^\star$ & iter.\\
\midrule
3 & 0.75 & -- & 0.0423 & 0.0380 & $0.1$ at $(0.55,0.5)$ & 2\\
9 & 0.75 & -- & 0.0424 & 0.0372 & $0.1$ at $(0.65,0.5)$ & 2\\
9 & 0.75 & 0.03 & 0.0424 & 0.0387 & $0.03$ at $(0.45,0.5),(0.55,0.5),(0.65,0.5)$; $0.01$ at $(0.55,0.7)$ & 2\\
9 & 1.5 & -- & 0.0805 & 0.0718 & $0.1$ at $(0.65,0.5)$ & 2\\
\bottomrule
\end{tabular}
\end{table}

\begin{table}[!t]
\centering
\caption{Scaling of the stress-chain reduction ($h=1/32$, $N=5$, $\rho=0.2m$), with the minimum-norm-point
algorithm for submodular minimization. Oracle calls are summed over observations, and ``bound'' is $N(2m-1)$.
Stress-corner evaluations are cached, so solve counts are distinct PDE solves. Times are total wall-clock seconds on
one CPU core.\label{tab:scaling}}
\small\setlength{\tabcolsep}{4.5pt}
\begin{tabular}{rrrcccrrcc}
\toprule
 & & & & oracle & & distinct & time & break- & \\
$m$ & $2^m$ & $N2^m$ & $|\mathcal C_n|$ & calls & bound & PDE solves & (s) & points & $|\operatorname{supp}P^\star|$\\
\midrule
8 & 256 & 1{,}280 & 5--6 & 43 & 75 & 178 & 2.2 & 24 & 6\\
12 & 4{,}096 & 20{,}480 & 6--8 & 55 & 115 & 724 & 7.8 & 30 & 6\\
16 & 65{,}536 & 327{,}680 & 8--10 & 75 & 155 & 2{,}186 & 23.7 & 40 & 6\\
20 & 1{,}048{,}576 & 5{,}242{,}880 & 10--11 & 93 & 195 & 5{,}455 & 59.0 & 49 & 6\\
\bottomrule
\end{tabular}
\end{table}

\begin{table}[!t]
\centering
\caption{Mesh refinement for the experiment of \cref{tab:dro,tab:solver} ($m=4$, $N=8$, same data, three
sites). On every grid, the stress chains coincide with those in \cref{tab:dro}, and $\Lambda$ has $19$ breakpoints,
whose positions vary with $h$. On the three grids where it was computed (not $h=1/256$), the optimal design assigns
all capacity to the station at $(0.55,0.5)$ in two cutting-plane iterations. The last row is a first-order
extrapolation to $h=0$ (see text).\label{tab:mesh}}
\small
\begin{tabular}{lcccccc}
\toprule
$h$ & $J_h(0.75)$ & $\lambda^\star_h$ & $J_h(1.5)$ & $\lambda^\star_h$ & $R_{0.75,h}(y^\star)$ & $R_{1.5,h}(y^\star)$\\
\midrule
$1/32$  & $0.04978$ & $0.06286$ & $0.09365$ & $0.05127$ & $0.03545$ & $0.06803$\\
$1/64$  & $0.05292$ & $0.06664$ & $0.09936$ & $0.05425$ & $0.03801$ & $0.07316$\\
$1/128$ & $0.05529$ & $0.06957$ & $0.10372$ & $0.05654$ & $0.03985$ & $0.07687$\\
$1/256$ & $0.05548$ & $0.06978$ & $0.10407$ & $0.05673$ & -- & --\\
\midrule
$h\to0$ (extrapolated) & $0.0557$ & & $0.1044$ & & & \\
\bottomrule
\end{tabular}
\end{table}

\paragraph{Stress chains and the exact robust value} Since $m=4$, we implement the oracle of \cref{alg:chain} by
enumeration. The chains contain two to four sets, below the limit of five, and the number of oracle calls equals the
bound $2|\mathcal C_n|-3$ for every observation (\cref{tab:dro}). For every radius, only $27$ chain corners are needed,
compared with $N2^m=128$ corners. At $401$ values of $\lambda$, the largest maximizers computed by enumeration are
nested and belong to the computed chains. Ties are well separated: exactly two sets tie at each breakpoint, every
other set is lower by at least $2.2\times10^{-4}$, and the chains do not change when the tolerance varies from
$10^{-14}$ to $10^{-4}$. The construction of \cref{thm:exact} agrees with a linear program to $2\times10^{-16}$ at $24$
radii, including $\rho=0$, every kink of $J_h$ (at each minimizer $\lambda^\star$), and radii beyond the total
budget; a variant with observations on the upper boundary, for which $K_n\neq\emptyset$ and several sets tie, passes
the same checks. Linear programs over distributions supported on (a) the $128$ corners and (b) these corners
together with a $5^4$-point grid of $\Xi$ agree with~\eqref{eq:exact} to $5.6\times10^{-17}$, so additional support
points never increase the worst-case value. (The $128$ corners contain $121$ distinct points, since the eight full
corners all equal $\xib$, which is also a grid point; hence $745$ distinct PDE solves.) Sources $2$ and $3$, which
are aligned with the protected zone, never leave the stressed set before sources $1$ and $4$, and in seven of the
eight observations they leave last (\cref{fig:dro}, left). At $\rho=0.75$, $\lambda^\star=0.0666$, and the worst-case
distribution moves observations $2$ and $3$ to their $\{2,3\}$-corners and splits observation $4$ with $p_4=0.548$,
as in~\eqref{eq:worstcase}. The worst-case expected damage is about $26$ times the empirical mean $0.0020$ at
$\rho=0.75$ and $49$ times at $\rho=1.5$.

\paragraph{Remediation design} Since $F$ is nonincreasing in $y$, an optimal design under the budget
$\sum_jy_j\le0.1$ can be chosen to use the full capacity, so we impose $\sum_jy_j=0.1$. For three sites, evaluating
all $45$ allocations on a simplex grid of step $0.0125$ ($5{,}760$ PDE solves) gives a smallest concentration of
$1.26\times10^{-4}$. The uniform test after \cref{prop:positivity} gives $u^n+0.1\min_jv^n_j\ge1.26\times10^{-4}$ for
all $n$, so $Y_+=Y$: every state that determines $R_{\rho,h}$ is nonnegative for every three-site design, here and on
the grids $h=1/32$ and $h=1/128$ of \cref{tab:mesh}. For the nine-site layout, the uniform test fails (minimum
$-1.4\times10^{-2}$), but the computed optimal designs satisfy the constraints defining $Y_+$ with a margin of
$1.26\times10^{-4}$, so they also solve the problem restricted to $Y_+$. By contrast, at $\xi=0$, placing all capacity
at one station gives concentrations as low as $-0.18$. All $84$ second differences of $R_{0.75,h}$ and $R_{1.5,h}$
along the $21$ grid lines with at least three points are positive, in agreement with \cref{prop:design}(c).

We then minimize $R_{\rho,h}$ over $Y=\{y\ge0,\ \sum_jy_j=0.1\}$, with $y_j\le\bar y$ when each site has a capacity
limit, by Kelley's cutting-plane method \cite{Kelley1960}. Given $R_{\rho,h}(y^i)$ and subgradients
$\sigma^i\in\partial R_{\rho,h}(y^i)$, $i=1,\dots,r$, the next iterate solves the master problem
\[
  \min_{y\in Y,\ \theta\in\R}\ \theta\quad\text{subject to}\quad\theta\ge R_{\rho,h}(y^i)+\ip{\sigma^i}{y-y^i},\qquad i=1,\dots,r.
\]
Its value is a lower bound (LB), and $\min_iR_{\rho,h}(y^i)$ is an upper bound (UB). The subgradients come
from~\eqref{eq:robustgrad}, and the adjoint gradients agree with finite differences to $5\times10^{-8}$. For three
sites, the method finds the grid-search optimum ($R_{0.75,h}=0.03801$), with all capacity at the central station
(\cref{tab:solver}). In every run, the second iterate $y^\star$ satisfies $0\in\sigma^\star+N_Y(y^\star)$ for its
subgradient $\sigma^\star$ up to rounding, and the master gap is at the level of rounding error
($|\mathrm{UB}-\mathrm{LB}|\le1.4\times10^{-17}$, relative $2\times10^{-16}$); overall accuracy is limited by the PDE and
adjoint solves (relative Newton residual $10^{-12}$). In these examples, the optimal designs are vertices of $Y$, and
the method terminates after two iterations. With nine candidate sites, the method moves capacity downstream along the
plume axis and reduces the robust risk to $0.0372$; per-site limits spread the allocation along that axis. The
nominal risk of the three-site optimum, $2.5\times10^{-4}$, is more than two orders of magnitude smaller than
$R_{0.75,h}$.

\paragraph{Scaling} To test larger problems, we repeat the experiment on a coarser grid ($h=1/32$) with
$m\in\{8,12,16,20\}$ sources in two staggered columns, $N=5$ observations, source amplitudes $100/m$, and $\rho=0.2m$,
so that the total source strength and the budget per coordinate ($\rho/m=0.2$) remain comparable. For submodular
minimization we use the minimum-norm-point method \cite{Wolfe1976,FujishigeIsotani2011}, which is efficient in
practice; the bound of \cref{prop:algorithm} also holds for strongly polynomial algorithms
\cite{Schrijver2000,IwataFleischerFujishige2001,Orlin2009}. The method stops when Wolfe's optimality test holds with
relative tolerance $10^{-12}$. The smallest and largest minimizers are then read off from the minimum-norm base $x$ as
$\{k:x_k<-\tau\}$ and $\{k:x_k\le\tau\}$ with $\tau=10^{-9}\max(1,\|x\|_\infty)$, and $f(S)-\sum_k\min(x_k,0)$
certifies optimality for the computed values of $F_h$, which carry the Newton error (largest gap
$8.5\times10^{-15}$). \Cref{alg:chain} accepts a new chain set only if it improves $\Psi_n$ by more than
$10^{-12}\max(1,|\Psi_n|)$. In \cref{tab:scaling}, the chains have about $m/2$ sets, well below $m+1$ (at most $11$
versus $21$ at $m=20$), and the number of oracle calls equals $2|\mathcal C_n|-3$ for every observation. The number of
distinct PDE solves grows from $178$ at $m=8$ to $5{,}455$ at $m=20$, whereas checking every corner at $m=20$ would
require $5{,}242{,}880$ solves, about $960$ times as many. For $m\le12$, direct enumeration confirms every chain, and
hence the tie handling, and reproduces $J_h(\rho)$ to $1.5\times10^{-15}$. In all tests, the worst-case distribution
has $N+1$ support points and splits one observation. For $m=4$, the grid $h=1/32$ already gives the same chains as
the finer grids (\cref{tab:mesh}).

\paragraph{Mesh refinement} We repeat the experiment of \cref{tab:dro,tab:solver} with the same data on grids with
$h=1/32$, $1/64$, $1/128$, and $1/256$ (\cref{tab:mesh}). The combinatorial structure does not change: on all four
grids, the chains of all eight observations coincide and $\Lambda$ has $19$ breakpoints (their positions move with
$h$; $\lambda^\star_h$ for $\rho=0.75$ increases from $0.0629$ to $0.0698$), the linear-program check
of~\eqref{eq:exact} agrees to $10^{-17}$, and, on the three grids where it was computed, the optimal design assigns
all capacity to the central station for both radii. All concentrations at the tested stress scenarios are positive.
The robust values increase under refinement, but the successive differences of $J_h(0.75)$, namely
$3.1\times10^{-3}$, $2.4\times10^{-3}$, and $1.8\times10^{-4}$, do not decrease regularly. Two first-order effects
contribute: upwinding adds the numerical diffusion $\beta_1h/2$, which is $31\%$ of $D$ at $h=1/64$, and $Z$ is
represented by the grid nodes it contains, so its effective area exceeds $|Z|=0.06$ by $2.5\%$, $0.5\%$, $3.1\%$, and
$1.8\%$ on the four grids. If $Z$ is represented instead by the exact fraction of each grid cell that it covers, the
successive differences decrease by factors of $1.91$ and $1.99$ for $\rho=0.75$ ($1.91$ and $2.00$ for $\rho=1.5$), as
expected for a first-order scheme, and the chains coincide with those of \cref{tab:dro} for $h\le1/64$. First-order
extrapolation from the two finest grids gives $J(0.75)\approx0.0557$ and $J(1.5)\approx0.1044$ for both
representations of $Z$. Thus the values at $h=1/64$ are about $5\%$ below the extrapolated estimates for the
continuous model, while the chains, the number of breakpoints, and, where computed, the optimal design are the same
on all grids tested.

\FloatBarrier
\section{The fractional $p$-Laplacian: a counterexample}\label{sec:plaplace}

We next examine whether \cref{thm:damage} remains valid when the fractional Laplacian $(-\Delta)^s$ is replaced by the
fractional $p$-Laplacian $(-\Delta_p)^s$. The following counterexample shows that it fails in general for $p>2$, even
without a reaction term ($g\equiv0$).

\begin{proposition}\label[proposition]{prop:plaplace}
Let $\Omega\subset\R^d$ be bounded and open, let $p\in(1,\infty)$ and $s\in(0,1]$, where $s=1$ means the local
$p$-Laplacian, and let $f\in L^\infty(\Omega)$ with $f\ge0$ and $f\ne0$. For $t\ge0$, let $u_t$ be the weak solution of
$(-\Delta_p)^su=tf$ in $\Omega$, with $u=0$ in $\R^d\setminus\Omega$ if $0<s<1$ and $u=0$ on $\partial\Omega$ if
$s=1$. Then $u_t=t^{1/(p-1)}u_1$, with $u_1\ge0$ and $u_1\ne0$. Consequently, for $f_1=f_2=f$ and $\Phi(r)=r$, we
have $F(\xi)=(\xi_1+\xi_2)^{1/(p-1)}\int_\Omega u_1$, and, for $\xi\in\R^2_+$ and $\delta>0$,
\[
  F(\xi+\delta e_1+\delta e_2)+F(\xi)-F(\xi+\delta e_1)-F(\xi+\delta e_2)\ \begin{cases}<0,&p>2,\\=0,&p=2,\\>0,&1<p<2.\end{cases}
\]
\end{proposition}
\begin{proof}
Existence and uniqueness follow from the strict monotonicity and coercivity of $(-\Delta_p)^s$ on the energy
space \cite{Lindqvist2017,IannizzottoMosconiSquassina2016}. Testing the equation with the negative part $u_1^-$ and
using $|a-b|^{p-2}(a-b)(a^--b^-)\le-|a^--b^-|^p$ (or $|\nabla u|^{p-2}\nabla u\cdot\nabla u^-=-|\nabla u^-|^p$ when
$s=1$) gives $u_1\ge0$. Since $f$ is not identically zero, neither is $u_1$. The operator is homogeneous of degree $p-1$, so
$t^{1/(p-1)}u_1$ solves the equation with load $tf$, and uniqueness gives $u_t=t^{1/(p-1)}u_1$. The second
difference equals $c[G(T+2\delta)+G(T)-2G(T+\delta)]$, where $c=\int u_1>0$, $T=\xi_1+\xi_2$, and
$G(\tau)=\tau^{1/(p-1)}$. The claim follows because $G$ is strictly concave when $p>2$ and strictly convex when
$1<p<2$.
\end{proof}

\begin{remark}
The proof of \cref{thm:order} does not extend directly to the $p$-Laplacian. The equations for differences of solutions have coefficients that depend on those solutions, so the two differences $d$ and $e$ used
in~\eqref{eq:de} do not satisfy equations with the same principal operator. \Cref{prop:plaplace} also shows that a
result like \cref{thm:damage}, for all nondecreasing convex damage functions, can hold only when $1<p\le2$. Any
positive result for $p>2$ would need additional assumptions on the damage or the sources.
\end{remark}

\section{Conclusion and open questions}\label{sec:open}
We have shown that comparison properties of elliptic PDEs can simplify Wasserstein distributionally robust
optimization. When the reaction is nondecreasing and concave, the load-to-damage function is nondecreasing, convex,
and supermodular. With a weighted $\ell_1$ transport cost, this structure reduces the worst-case optimization problem
to nested stress corners. The exact robust value can be computed from $\lambda=0$ and at most $Nm$ positive
breakpoints. Moreover, for every Wasserstein radius, a worst-case distribution can be chosen from a fixed set of at
most $N(m+1)$ chain corners.

Under nonnegative loading and $g(0)=0$, concavity of the reaction is also necessary for these properties to hold for
all admissible operators and damage functions. For a fixed damage function, \cref{cor:tradeoff} gives a weaker
necessary condition in the large-domain limit. We applied the results to convex remediation-capacity planning with
Monod biodegradation. In synthetic experiments, the stress-chain method required far fewer PDE solves than full
corner enumeration. The stress chains remained unchanged across the tested grids, and the optimal design agreed on
the grids for which it was computed. We also gave a counterexample showing that supermodularity can fail for the
fractional $p$-Laplacian when $p>2$.

Two questions remain open. First, for a convex reaction, is condition~\eqref{eq:tradeoff}, or an analogous condition
for nonlocal operators, sufficient for supermodularity of $F$ on general domains? Second, when $1<p<2$ and
$g\equiv0$, is the solution map supermodular for nonproportional sources? It is also of interest to determine how a
concave reaction changes this second question, since the homogeneity used in \cref{prop:plaplace} is then lost.
\appendix
\section{Verification of \texorpdfstring{\cref{ass:H}}{the operator assumption} for the examples}\label[appendix]{app:examples}
\paragraph{\Cref{ex:local}} The space $H^1_0(\Omega)$ is continuously and densely embedded in $L^2(\Omega)$. For
$u\in H^1_0$, integration by parts gives
$\int(\beta\cdot\nabla u)u=\frac12\int\beta\cdot\nabla(u^2)=-\frac12\int(\operatorname{div}\beta)u^2$, so
$\E(u,u)\ge\lambda_A\|\nabla u\|_{L^2}^2+\int(a-\frac12\operatorname{div}\beta)u^2\ge\lambda_A\|\nabla u\|_{L^2}^2$, where $\lambda_A>0$ is the ellipticity constant of $A$. Poincar\'e's inequality now gives coercivity. For $v\in H^1_0$, we have $v^\pm\in H^1_0$, with
$\nabla v^+=\mathbf 1_{\{v>0\}}\nabla v$ and $\nabla v^-=-\mathbf 1_{\{v<0\}}\nabla v$
\cite[Lemma~7.6]{GilbargTrudinger2001}. Each integrand in $\E(v^+,v^-)$ contains one factor from
$\{\nabla v^+,v^+\}$ and one from $\{\nabla v^-,v^-\}$, and these have disjoint supports almost everywhere.
Therefore, $\E(v^+,v^-)=0$.

\paragraph{\Cref{ex:fractional}} The space $X^s_0(\Omega)$ is a Hilbert space with the Gagliardo seminorm; it is
continuously embedded in $L^2(\Omega)$ and contains $C_c^\infty(\Omega)$ \cite{ServadeiValdinoci2013}. The fractional
Poincar\'e inequality and $a\ge0$ give coercivity. Since $|v^+(x)-v^+(y)|\le|v(x)-v(y)|$, taking positive parts does
not increase the Gagliardo seminorm, so $v^+\in X^s_0$. Pointwise,
$(v^+(x)-v^+(y))(v^-(x)-v^-(y))=-v^+(x)v^-(y)-v^+(y)v^-(x)\le0$ and $v^+v^-=0$; hence $\E(v^+,v^-)\le0$.

\paragraph{\Cref{ex:graph}} Since the graph is finite, $V$ is finite-dimensional. By the absorption assumption, the
energy $\E(v,v)=\frac12\sum\omega_{xy}(v_x-v_y)^2+\sum\kappa_xv_x^2$ vanishes only when $v=0$, so $\E$ is coercive. The pointwise argument from \cref{ex:fractional} also gives $\E(v^+,v^-)\le0$.

\paragraph{\Cref{ex:discrete}} We have $\E(v^+,v^-)=\sum_{i,j}v^-_iK_{ij}v^+_j=\sum_{i\ne j}K_{ij}v^-_iv^+_j\le0$, because
the diagonal terms vanish ($v^-_iv^+_i=0$) and $K_{ij}\le0$ for $i\ne j$. This gives \ref{H3}. Coercivity follows
from $v^\top Kv=\frac12v^\top(K+K^\top)v$, since the symmetric part of $K$ is positive definite. For the scheme of \cref{sec:numerics}, let $K_h^{\rm FD}=D\,L_h+B_h+\gamma I$ be the finite-difference matrix of the
strong form, where $L_h$ is the five-point approximation of $-\Delta$ (a symmetric positive definite Z-matrix) and
$B_h=\frac b h\,(I-S)\otimes I$ is upwind advection for $\beta=(b,0)$, $b\ge0$, with $S$ the lower shift. Then
$K_h^{\rm FD}$ is a Z-matrix, and $B_h+B_h^\top=\frac bh\,(2I-S-S^\top)\otimes I$ is positive semidefinite, so
$K_h^{\rm FD}+(K_h^{\rm FD})^\top$ is positive definite. With the lumped weights $\mu_i=h^2$, that is, $M_h=h^2I$, the
bilinear form is $\E_h(u,\varphi)=\varphi^\top M_hK_h^{\rm FD}u$, so $K=h^2K_h^{\rm FD}$ in \cref{ex:discrete}, and
multiplication by $h^2>0$ preserves both properties. The reaction and load terms are $\varphi^\top M_hg(u)$ and
$\varphi^\top M_hf$, so~\eqref{eq:state} is equivalent to $K_h^{\rm FD}u+g(u)=f$ at the grid nodes, which is the system
solved in \cref{sec:numerics}. For monotone
discretizations of the fractional Laplacian \cite{HuangOberman2014}, the same verification applies once the
matrix is shown to be a Z-matrix with positive definite symmetric part.

\section{Observations used in \texorpdfstring{\cref{sec:numerics}}{Section 7}}\label[appendix]{app:data}
\Cref{tab:obs} lists the observations $\xih_n$ used in the experiments with $m=4$ (\cref{tab:dro,tab:solver,tab:mesh}),
rounded to four digits; the code uses the unrounded values generated with seed $7$.

\begin{table}[h]
\centering
\caption{Observations $\xih_n$ for $m=4$, $N=8$, and $\Xi=[0,2]^4$.\label{tab:obs}}
\small
\begin{tabular}{ccccc}
\toprule
$n$ & $\xih_{n1}$ & $\xih_{n2}$ & $\xih_{n3}$ & $\xih_{n4}$\\
\midrule
$1$ & $0.9004$ & $0.9992$ & $0.8177$ & $0.6590$\\
$2$ & $0.7676$ & $0.6361$ & $0.9191$ & $1.4387$\\
$3$ & $0.7576$ & $0.7243$ & $1.0683$ & $1.0197$\\
$4$ & $0.9338$ & $0.6498$ & $0.8908$ & $1.1480$\\
$5$ & $0.5622$ & $0.7668$ & $0.4626$ & $0.5731$\\
$6$ & $0.4724$ & $0.8289$ & $0.5775$ & $0.9896$\\
$7$ & $0.9508$ & $0.8430$ & $0.3730$ & $0.7453$\\
$8$ & $0.8849$ & $0.9364$ & $0.5268$ & $0.7614$\\
\bottomrule
\end{tabular}
\end{table}

\section*{Declarations}

\smallskip
\noindent\textbf{Funding.}
This work was supported by the Department of Science and Technology (DST),
Government of India, through the INSPIRE Faculty Award, grant
no.\ IFA21-MA 158.

\medskip
\noindent\textbf{Data availability.}
The code and data used to reproduce the numerical experiments in \cref{sec:numerics} are
available from the corresponding author upon reasonable request.

\medskip
\noindent\textbf{Conflict of interest.}
The authors have no competing interests to declare that are relevant to the
content of this article.

\raggedbottom

\end{document}